\documentclass{amsproc}
\usepackage{amsfonts}
\usepackage{amssymb}
\usepackage{amscd}
\usepackage{verbatim}
\usepackage{graphicx}
\usepackage{subfigure}
\usepackage{fancyhdr}
\usepackage[T1]{fontenc}
\usepackage[latin1]{inputenc}
\usepackage[english]{babel}
\usepackage{hyperref}
\usepackage[active]{srcltx}
\input xy
\xyoption{all}

\theoremstyle{plain}

\makeatletter

\@addtoreset{equation}{section} \makeatother

\newcommand{\fait}[3]{\begin{#1}\label{#2}{#3}\end{#1}}
\newcommand{\maps}[2]{$#1 \mapsto #2$}
\newcommand{\mor}[3]{$#1:#2\rightarrow #3$}

\newcommand{\iso}[3]{$\xymatrix@1@C=15pt{#3: #1\ar[r]^-{\simeq}& #2}$}

\newcommand{\rad}[1]{{\rm rad}#1}
\newcommand{\Hom}[1]{{\rm Hom}#1}
\newcommand{\Ext}[1]{{\rm Ext}#1}

\newcommand{\Ker}[1]{{\rm Ker}\ #1}

\newcommand{\SH}[1]{{\rm SH}#1}

\begin{document}

\title{(Co)Homology theories for oriented algebras}

\author[J.~C.~Bustamante]{Juan Carlos Bustamante}
\address{J.~C.~Bustamante; Departamento de Matem\'aticas, Universidad San Francisco de Quito, Quito, Ecuador.}
\email{juanb@usfq.edu.ec}

\author[J.~Dionne]{Julie Dionne}
\address{J.~Dionne; D\'epartement de Math\'ematiques, Universit\'e
de Sherbrooke, Sherbrooke, J1K 2R1, Qu\'ebec, Canada.}
\email{julie.dionne2@usherbrooke.ca}

\author[D.~Smith]{David Smith \footnote{Corresponding author : D.~Smith}}
\address{D.~Smith; Department of Mathematics, Bishop's University, Sherbrooke, Québec,
Canada, J1M0C8.}
\email{unracinois@hotmail.com}



\begin{abstract} We study three different (co)homology theories
for a family of pullbacks of algebras that we call oriented. We
obtain a Mayer Vietoris long exact sequence of Hochschild and cyclic
homology and cohomology groups for these algebras. We give examples
showing that our sequence for Hochschild cohomology groups is
different from the known ones. In case the algebras are given by
quiver and relations, and that the simplicial homology and
cohomology groups are defined, we obtain a similar result in a
slightly wider context. Finally we also study the fundamental groups
of the bound quivers involved in the pullbacks.
\end{abstract}

\maketitle


%
%
\section*{Introduction} \label{Introduction}
Let $k$ be a field and $R$ be an associative $k$-algebra with unit.
Given an $R-R$-bimodule $_RM_R$, or equivalently a right (or left)
module over the enveloping algebra $R^e:=R\otimes_kR^{op}$, the
Hochschild homology and cohomology groups of $R$ with coefficients
in $M_{R^e}$ are the groups ${\rm H}_n(R,M)={\rm Tor}^{R^e}_n(M,R)$
and ${\rm H}^n(R,M)={\rm Ext}_{R^e}^n(R,M)$, respectively. In case
$_RM_R=_R\!\!R_R$ we simply write ${\rm H}_n(R)$, and ${\rm
H}^n(R)$.

The low dimensional cohomological groups ${\rm H}^i(R)$ (that is
$i\leq 2$) have clear classical interpretations in terms of the
center, the derivations or the extensions of $R$. Besides this,
these groups also play an important role in representation theory of
algebras: the rigidity properties of $R$ are closely related to
${\rm H}^2(R)$ and ${\rm H}^3(R)$, whereas ${\rm H}^1(R)$ gives
information about the simple connectedness of $R$, and there are
links between the vanishing of ${\rm H}^1(R)$ and ${\rm H}^2(R)$ and
the representation type of $R$. Unfortunately, the groups ${\rm
H}^i(R)$ are not easy to compute in general. There are some long
exact sequences allowing to perform computations in some cases, see
\cite{Ci00, MP00, GS02, GMS03, BG04}.

The main motivation for this work is to compute the Hochschild
cohomology groups of the pullback $R$ of two morphisms of
$k$-algebras $A_1 \rightarrow C \leftarrow A_2$.  However, as
Examples \ref{ex tournicoti} show, there is no easy immediate
relation between the Hochschild cohomology groups of the algebras
$R,A_1, A_2$ and $C$. Thus, we impose additional hypotheses on the
morphisms above, and obtain a Mayer-Vietoris-like long exact
sequence of Hochschild (co)homology groups of what we call {\bf
oriented algebras}.  While doing so,  one can see that using the
same techniques it is possible to obtain analogous results about
cyclic (co)homology as well as for  simplicial (co)homology groups
of algebras (if they are defined) and fundamental groups of the
involved bound quivers, if the algebras are given in this way. This
paper is organised as follows:

In Section \ref{Preliminaries} we recall some basic constructions
related to bound quivers, and algebras given by the latter. We then
define {oriented algebras}, give some examples, and make some
preliminary observations.

Section \ref{Hochschild (co)homology} is devoted to Hochschild
(co)homology. After recalling the basic definitions we establish the
existence of two long exact sequences of Hochschild homology and
cohomology groups, in Theorems \ref{thm Hochschild homology} and
\ref{thm Hochschild cohomology}.  Moreover, we show that the
morphisms involved in the cohomological sequence are compatible with
the Gerstenhaber algebra structure of ${\rm
H}^*(R):=\bigoplus_{i\geq 0} {\rm H}^i(R)$. The remaining part of
this section is devoted to specialize our results to the case where
the algebra $C$ is what we call a core algebra, and, on the other
hand, to  compare the long exact sequence of Theorem \ref{thm
Hochschild cohomology} with the known long exact sequences of
Hochschild cohomology groups of \cite{Ci00, MP00, GS02, GMS03,
BG04}.

In Section \ref{Cyclic (co)homology}, again, after recalling some
constructions, we establish the existence of long exact sequences of
cyclic homology and cohomology groups, in Theorems \ref{thm cyclic
homology} and \ref{thm cyclic cohomology}. Moreover, we show that
these sequences are compatible with Connes' long exact sequences
relating Hochschild and cyclic (co)homology groups.

Finally, in Section  \ref{Simplicial and fundamental} we turn our
attention to algebras given by bound quivers, and, more precisely to
algebras having semi-normed basis. The existence of such basis
allows to define the simplicial (co)homology groups of algebras,
following  \cite{MdlP99}. In the same flavor of the previous parts
of this paper, but with weaker hypotheses, we establish the
existence of long exact sequences of simplicial (co)homology groups,
in Propositions \ref{SimpHom} and \ref{SimpCohom}. We then consider
fundamental groups of bound quivers, and obtain an explicit formula
that allows to compute them  in a context which is slightly
different from that of the oriented pullbacks.

%
%
\section{Preliminaries}
    \label{Preliminaries}

\subsection{Algebras and quivers}
    \label{Algebras and quivers}
Throughout this paper, $k$ denotes a commutative field and all
algebras are associative $k$-algebras with unit. Tensor products are
taken over $k$ unless otherwise stated, so $\otimes=\otimes_k$. Let
$R^{op}$ be the opposite algebra of $R$, and denote by $R^e$ its
enveloping algebra $R\otimes R^{op}$ of $R$. We shall freely use the
fact that the category of $R$-bimodules is equivalent to the
category of right (or left) $R^e$-modules: given a $R$-bimodule $M$,
define the right $R^e$-structure by  $m(a\otimes b)=bma$ (and the
left $R^e$-structure by $(a\otimes b)m=amb$).

A (finite) \textbf{quiver} is a quadruple $Q=(Q_0, Q_1, s, t)$,
where $Q_0$ and $Q_1$ are two (finite) sets, respectively called the
\textbf{set of vertices} and the \textbf{set of arrows}, and $s$ and
$t$ are two maps from $Q_1$ to $Q_0$ associating to each arrow its
\textbf{source} and its \textbf{target}, respectively. A quiver $Q$
is \textbf{connected} if its underlying graph is connected.

A \textbf{path} on \textbf{length} $n$  is a sequence of $n$  arrows
$\omega = \alpha_1 \alpha_2 \cdots \alpha_n$ such that $t(\alpha_i)
= s(\alpha_{i+1})$ for each $i$, with $1 \leq i \leq n-1$.  The
\textbf{source} of such a path is $s(\alpha_1)$, and its
\textbf{target} is $t(\alpha_n)$. Also, to each vertex $x$ of $Q$,
we associate the \textbf{stationary path} $\varepsilon_x$, whose length is $0$, and whose source and target are $x$. An
\textbf{oriented cycle} is a path $\omega=\alpha_1 \alpha_2 \cdots \alpha_n$ of length at least one, such that
$s(\alpha_1)=t(\alpha_n)$.

The \textbf{path algebra} $kQ$ is the $k$-vector space having as
basis the set of all paths in $Q$ endowed with the following
multiplication: let $\omega_1$ and $\omega_2$ be two basis elements,
then $\omega_1 \omega_2$ is their composition if $t(\omega_1) =
s(\omega_2)$, and 0 otherwise. A \textbf{relation} from $x$ to $y$,
with $x, y \in Q_0$, is a linear combination $\rho = \sum_{i=1}^r
\lambda_i \omega_i$ where, for each $i$, $\omega_i$ is a path of
length at least two from $x$ to $y$ in $Q$ and $\lambda_i\in
k\setminus \{0\}$. Let $F$ be the two-sided ideal of $kQ$ generated
by all the arrows of $Q$. An ideal $I$ of $kQ$ is
\textbf{admissible} if there exists $m \geq 2$ such that $F^m
\subseteq I \subseteq F^2$. In this case, the pair $(Q, I)$ is
called a \textbf{bound quiver}. The algebra $kQ/I$ is basic,
connected (if so is $Q$), and finite dimensional (if $Q$ is finite).
By abuse of notation, we simply denote the primitive pairwise
orthogonal idempotents $\varepsilon_x +I $ of $kQ/I$ by
$\varepsilon_x$, so that $1_{kQ/I}= \sum_{x\in Q_0} \varepsilon_x$.

Conversely, for every finite dimensional, connected and basic
algebra $A$ over an algebraically closed field $k$, there exists a
unique connected quiver $Q$ and surjective maps \mor{\nu}{kQ}{A}
with $I= \Ker \nu$ admissible (see \cite{BG82}). The pair $(Q, I)$
is called a \textbf{presentation} of $A$. Following \cite{BG82}, any
basic algebra $A=kQ/I$ can equivalently be regarded as a locally
bounded $k$-category  having as object class the set $A_0 = Q_0$ and
as morphism set from $x$ to $y$ the $k$-vector space
$A(x,y)=\varepsilon_xA\varepsilon_y$. In what follows, we use both
terminologies, and pass from one to the other freely. A subcategory
$B$ of $A$ is called \textbf{full} if, for all $x,y \in B_0$, we
have $B(x,y) = A(x,y)$, and \textbf{convex} if any vertex  which
lies on a path in $Q$ having its source and target in $B_0$ also
belongs to $B_0$. Note that convexity does not imply fullness.

\subsection{Pullbacks} \label{pullbacks}
Let \mor{f_1}{A_1}{C} and \mor{f_2}{A_2}{C} be two surjective
morphisms of $k$-algebras. Then the set $R=\{(a_1,a_2)\in A_1\times
A_2| \  f_1(a_1)=f_2(a_2)\}$, endowed with the natural operations
becomes a $k$-algebra, which is in fact the pullback of $f_1$ and
$f_2$. As Examples \ref{ex tournicoti} below show, one cannot hope
to have a direct relation between the Hochschild cohomology groups
of $R$, $A_1$, $A_2$, and $C$ in general.

In \cite{IPTZ87}, a special kind of pullbacks has been considered:
assume that $1_{A_1}=u'_1+u''_1,\ 1_{A_2}=u'_2+u''_2$ and that there
exists an isomorphism $A_1/<u'_1> \simeq A_2/<u'_2>$ sending the
image of $u''_1$  under the projection, to the image of $u''_2$. Let
$C$ be this common quotient,  $f_1$ and $f_2$ the projections, and,
as above, $R$ be the pullback of $f_1$ and $f_2$. In $R$, let
$e'_1=(u'_1,0)$,  $e'_2=(0, u'_2)$, and $e=(u''_1, u''_2)$. A direct
computation shows that $e'_1 R e'_2 = e'_2 R e'_1=0$ and that
$1_R=e'_1+e+e'_2$. In \cite{IPTZ87} one can find the detailed
construction of the quiver of $R$ from those of $A_1, A_2$ and $C$.

\fait{examples}{ex tournicoti}{\emph{%
\begin{enumerate}
\item [(a)] Let $A_1$ and $A_2=A_1^{op}$  be the hereditary algebras
given by the quivers $\SelectTips{eu}{10}\xymatrix{1\ar[r]^{\alpha}
&2}$ and $\SelectTips{eu}{10}\xymatrix{1 &2\ar[l]_{\beta}}$,
respectively. Moreover, let $C$ be the semi-simple algebra generated
by the vertices $1$ and $2$, and, for $i\in\{1,2\}$ let
\mor{f_i}{A_i}{C} be the projections.  The pullback $R$ is then the
bound quiver algebra given by the quiver
$\SelectTips{eu}{10}\xymatrix{1\ar@<.5ex>[r]^{\alpha}
&2\ar@<.5ex>[l]^{\beta}}$ with relations $\alpha \beta =0= \beta
\alpha$.  For each $n> 0$, we have ${\rm H}^n(R) \neq 0$ by
\cite[(p.~96)]{Ci98}, while ${\rm H}^n(A_1)= {\rm H}^n(A_2)= {\rm
H}^n(C)=0$.
\item[(b)] Let $R$ be the algebra given by the quiver
$$\SelectTips{eu}{10}\xymatrix@R=10pt@C=10pt{& 2\ar[dr]^{\beta} \\
1\ar[ur]^{\alpha}&& 3\ar[dl]^{\gamma}\\
& 4\ar[ul]^{\delta}}$$ bounded by the ideal generated by all paths
of length two. Let $A_1$ and $A_2$ be full subcategories of $R$
generated by the sets of vertices $\{1,\ 2,\ 4\}$, and $\{2,\ 3,\
4\}$ respectively. Then, $R$ is the pullback  of the projections of
$A_1$ and $A_2$ over the semisimple algebra $C$ generated by the
vertices $1$ and $3$. By setting $e=\varepsilon_2+\varepsilon_4$,
$e'_1=\varepsilon_1$ and $e'_2=\varepsilon_3$, it is easily seen
that $1_R=e'_1+e+e'_2$ and $e'_1Re'_2=0=e'_2Re'_1$. %
\newline
Again, as in the previous example, we get from \cite[(p.~96)]{Ci98} that ${\rm H}^n(R) \neq 0$ for infinitely many $n>0$, while ${\rm H}^n(A_1)= {\rm H}^n(A_2)= {\rm H}^n(C)=0$ for each $n>0$. In fact what  allows to construct nontrivial elements in ${\rm H}^{4n}(R)$ is the fact that $eRe'_i \not= 0$ and $e'_iRe \not= 0$ for $i\in \{1,2\}$.
\end{enumerate}}}

The last example motivates the following definition.

%

\fait{definition}{defn oriented}{Let $R$ be an algebra. Assume that
there exists a decomposition $1_R=e'_1 + e + e'_2$ of the unit of
$R$, where $e, e'_1, e'_2$ are orthogonal idempotents of $R$.  Then
$R$ is said to be \emph{\textbf{oriented}} by $e, e'_1, e'_2$ if
\begin{eqnarray} \label{case 0}
e'_1Re'_2 = e'_2Re'_1 = 0
\end{eqnarray}
and one of the following conditions is satisfied :
\begin{eqnarray}
eRe'_1 = eRe'_2 = 0; \label{case 1} \\
e'_1Re = e'_2Re =0; \label{case 2} \\
e'_1Re = eRe'_2 =0.\label{case 3}
\end{eqnarray}
In this case, we set $e_i=e'_i + e, A'_i=e'_iRe'_i$, $A_i=e_iRe_i$
for $i=1,2$ and $C=eRe$. Finally, we set $E=ke'_1+ke+ke'_2\subseteq
R$, $E_1=ke'_1+ke\subseteq A_1$, $E_2=ke'_2+ke\subseteq A_2$ and
$E_C=ke\subseteq C$.}

Clearly, $R$, $C$, $A_1$ and $A_2$ are $E$-bimodules, and $C$, $A_1$
and $A_2$ are also bimodules over the semisimple algebras $E$, $E_1$
and $E_2$ respectively. As we shall see, these semisimple algebras
will play a crucial role in the sequel, as nice substitutes for $k$.

Outside the context of bound quivers, an algebra $\Lambda$ is said
to be {\bf convex} in bigger algebra $\Gamma$ ($\Lambda \subseteq
\Gamma$, but not necessarily a subalgebra) if whenever $\lambda_1,
\lambda_2, \ldots, \lambda_n$ are elements in $\Gamma$ such that
$1_\Lambda \lambda_1 \lambda_2 \cdots \lambda_n 1_\Lambda \neq 0$,
then each of the $\lambda_i$'s belong to $\Lambda$.

Also observe that conditions (\ref{case 0})-(\ref{case 3}) are
chosen in such a way that if $R$ is oriented by $e, e'_1, e'_2$,
then the algebras $A_i$, $A'_i$ and $C$ are convex in $R$, for
$i=1,2$.

\fait{examples}{rem oriented}{\emph{
\begin{enumerate}
\item[(a)]  If $e=0$ but $e'_1\neq 0$ and $e'_2\neq 0$,
then condition (\ref{case 0}) simply means that $R\simeq A'_1 \times
A'_2=A_1 \times A_2$.
\item[(b)]  Let $\Gamma$, $\Lambda$ be two $k$-algebras and
$_\Gamma M_\Lambda$ be a $\Gamma - \Lambda$-bimodule. We can then
construct the {\bf triangular matrix ring}
$R=\left(\begin{array}{cc} \Gamma & 0 \\ _\Lambda M_\Gamma & \Lambda
\end{array}\right)$.
By taking $e'_1=\left(\begin{array}{cc}1_\Gamma&0\\ 0&0 \end{array} \right),\ \ e'_2=0$
and $e=\left(\begin{array}{cc}0 & 0\\0&1_\Lambda \end{array} \right) $
one can see that $R$ is oriented. In case $\Lambda$ is a
division ring, the $R$ is called a {\bf one-point extension}
of $\Gamma$ by $M$. See Remark \ref{rem pullbacks}, (a), case 1.3.
\end{enumerate}}}

%

\fait{remark}{rem pullbacks}{\emph{%
\begin{enumerate}%
\item [(a)] If $R$ is oriented by $e, e'_1$ and $e'_2$,  using the decomposition
$1_R=e'_1 + e + e'_2$, the relation $R=1_R\cdot R\cdot 1_R$, and the
isomorphism of $k$-algebras $R\simeq {\rm End}_R(R)$, one can write
$R$ in a matrix form. The next figure shows the possible matrix
configurations, together with their corresponding quiver
configurations in case $R$ is a bound quiver algebra. Note that the
only possible nonzero paths in each case are indicated by the
depicted arrows.%
$$
{\setlength\arraycolsep{2pt}
\begin{array}{ccc}
  \left( \begin{array}{ccc} A'_1 & e'_1 R e & 0 \\
    0 & C &  0 \\
    0 & e'_2R e & A'_2
\end{array} \right) & \quad \quad \left( \begin{array}{ccc} A'_1 & 0  & 0 \\
    e R e'_1 & C & e R e'_2 \\
    0 & 0 & A'_2
\end{array} \right) \quad \quad &  \left( \begin{array}{ccc} A'_1 & 0 & 0 \\
    eRe'_1 & C & 0 \\
    0 & e'_2Re & A'_2
\end{array} \right) \\ \ \\
    \begin{picture}(0,0)%
\includegraphics{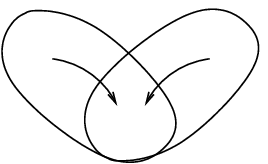}%
\end{picture}%
\setlength{\unitlength}{1740sp}%
\begingroup\makeatletter\ifx\SetFigFont\undefined%
\gdef\SetFigFont#1#2#3#4#5{%
  \reset@font\fontsize{#1}{#2pt}%
  \fontfamily{#3}\fontseries{#4}\fontshape{#5}%
  \selectfont}%
\fi\endgroup%
\begin{picture}(2836,1709)(1815,-6023)
\put(2026,-4741){\makebox(0,0)[lb]{\smash{{\SetFigFont{6}{7.2}{\rmdefault}{\mddefault}{\updefault}$A'_1$}}}}
\put(4051,-4741){\makebox(0,0)[lb]{\smash{{\SetFigFont{6}{7.2}{\rmdefault}{\mddefault}{\updefault}$A'_2$}}}}
\put(3151,-5641){\makebox(0,0)[lb]{\smash{{\SetFigFont{6}{7.2}{\rmdefault}{\mddefault}{\updefault}$C$}}}}
\end{picture}%

&
    \begin{picture}(0,0)%
\includegraphics{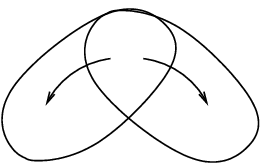}%
\end{picture}%
\setlength{\unitlength}{1740sp}%
\begingroup\makeatletter\ifx\SetFigFont\undefined%
\gdef\SetFigFont#1#2#3#4#5{%
  \reset@font\fontsize{#1}{#2pt}%
  \fontfamily{#3}\fontseries{#4}\fontshape{#5}%
  \selectfont}%
\fi\endgroup%
\begin{picture}(2836,1709)(1815,-5798)
\put(1981,-5371){\makebox(0,0)[lb]{\smash{{\SetFigFont{6}{7.2}{\rmdefault}{\mddefault}{\updefault}$A'_1$}}}}
\put(3916,-5371){\makebox(0,0)[lb]{\smash{{\SetFigFont{6}{7.2}{\rmdefault}{\mddefault}{\updefault}$A'_2$}}}}
\put(3106,-4561){\makebox(0,0)[lb]{\smash{{\SetFigFont{6}{7.2}{\rmdefault}{\mddefault}{\updefault}$C$}}}}
\end{picture}%

&
    \begin{picture}(0,0)%
\includegraphics{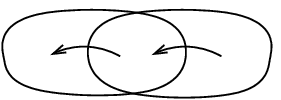}%
\end{picture}%
\setlength{\unitlength}{2030sp}%
\begingroup\makeatletter\ifx\SetFigFont\undefined%
\gdef\SetFigFont#1#2#3#4#5{%
  \reset@font\fontsize{#1}{#2pt}%
  \fontfamily{#3}\fontseries{#4}\fontshape{#5}%
  \selectfont}%
\fi\endgroup%
\begin{picture}(2555,863)(2751,-1593)
\put(2836,-1096){\makebox(0,0)[lb]{\smash{{\SetFigFont{7}{8.4}{\rmdefault}{\mddefault}{\updefault}$A'_1$}}}}
\put(3916,-1096){\makebox(0,0)[lb]{\smash{{\SetFigFont{7}{8.4}{\rmdefault}{\mddefault}{\updefault}$C$}}}}
\put(4771,-1096){\makebox(0,0)[lb]{\smash{{\SetFigFont{7}{8.4}{\rmdefault}{\mddefault}{\updefault}$A'_2$}}}}
\end{picture}%

\\ \ \\
  \mbox{case (\ref{case 1})} & \mbox{case (\ref{case 2})} & \mbox{case (\ref{case 3})}
\end{array}}
$$
\item  [(b)] Note that any matrix algebra of one of these forms is oriented
by the idempotents $e'_1=e_{11}, e=e_{22}, e'_2= e_{33}$, where
$e_{ii}$ is the $3\times 3$-matrix having a $1$ in position $(i,i)$
and $0$ elsewhere. Moreover, observe that in this case $R$ is the
pullback of the projections \mor{\pi_i}{A_i}{C} :
\maps{e_ire_i}{ere}, for $i=1,2$. Therefore, the oriented algebras
(with $e$, $e'_1$ and $e'_2$ nonzero) generalize the "Nakayama
oriented pullbacks" introduced in \cite{Le04}. However, not all
pullbacks are oriented, as Examples \ref{ex tournicoti} show.
\end{enumerate}}}

\subsection{Some immediate observations} \label{Immediate observations}
In this section, we establish some immediate results that will be
needed in the sequel. From now on we adopt the notations and
terminology of Definition \ref{defn oriented}.

\fait{remark}{rem tenseurs}{\emph{%
Let $R$ be an algebra oriented by a set of idempotents $e, e'_1, e'_2$.%
\begin{enumerate}
\item[(a)] Let $\sigma_1, \sigma_2, \tau_1, \tau_2 \in \{e'_1, e, e'_2\}$.
Then, in $R\otimes_E R$, we have $$ \sigma_1R\tau_1\otimes_{E} \sigma_2R
\tau_2= \sigma_1 R \tau_1 \sigma_2  \otimes_E \sigma_2 R \tau_2 =
\sigma_1 R  \tau_1 \otimes_E \tau_1 \sigma_2 R \tau_2$$ and, since
the idempotents $e'_1$, $e$ and $e'_2$ are orthogonal, this vanishes
if $\tau_1\neq \sigma_2$.  The same argument shows that $R^{\otimes_E n}$
is generated by the elements of the form
$r_1\otimes r_2\otimes \cdots\otimes r_n$, with  $r_i\in \sigma_i R \tau_i$
with $\sigma_i, \tau_i\in \{e, e'_1, e'_2\}$  and $\tau_i=\sigma_{i+1}$
for each $i$. Note that we are strongly using the fact that the tensor
products are taken over $E$. The claim does not hold if tensors are
taken over $k$.%
\item[(b)] As noted before, $C$ is full and convex in $A_1$ and $A_2$,
which in turn are full and convex in $R$. Thus, we have fully
faithful embeddings ${\rm mod}C \subseteq {\rm mod}A_i \subseteq {\rm mod}R$.
Consequently,  $\displaystyle X \otimes_E Y \simeq X \otimes_{E_C} Y$
whenever $X$ and $Y$ are modules over $C$ (or $C^{op}$).%
\item[(c)] Let $X$ be an $R$-bimodule, or, equivalently an $R^e$-module.
Then there exists an epimorphism \mor{p}{(R^e)^n}{X}, so that
for $x\in X$ we have $x=p(r)$, hence
$$e' _1xe'_2=(e'_1\otimes e'_2)x=(e'_1\otimes e'_2)p(r)=p((e'_1\otimes e'_2)r)= p(e'_1 r e'_2)=p(0)=0$$
This shows that  $e'_1Xe'_2=0$. Similarly, one can show that $e'_2Xe'_1=0$. %
\end{enumerate}}}

We then get the following lemmata

\fait{lemma}{radical2}{If $R$ is oriented by a set of idempotents
$e, e'_1$ and $e'_2$, then
\begin{enumerate}
    \item[(a)] $eR^{{\otimes}_E n}e =(eRe)^{{\otimes}_E n}=C^{{\otimes}_E n}$;
    \item[(b)] $e_1R^{{\otimes}_E n}e_1 =(e_1Re_1)^{{\otimes}_E n}=A_1^{{\otimes}_E n}$;
    \item[(c)] $e_2R^{{\otimes}_E n}e_2 =(e_2Re_2)^{{\otimes}_E n}=A_2^{{\otimes}_E n}$;
\end{enumerate}}

\begin{proof}
We only prove (a) since the proofs of (b)
and (c) are similar.\\
(a). Clearly, $eRe$ is included in $R$, whence $(eRe)^{{\otimes}_E
n} \subseteq R^{{\otimes}_E n}$. Multiplying by $e$ on both sides
gives $(eRe)^{{\otimes}_E n} \subseteq eR^{{\otimes}_E n}e$.
On the other hand, let $\tau = e(r_1 \otimes r_2
\otimes\cdots\otimes r_n)e$ be a nonzero element of $eR^{{\otimes}_E
n}e$. It then follows from Remark \ref{rem tenseurs} (a) and the
convexity of $C$ in $R$ that $r_i\in eRe$, for each $i$. Thus
$\tau\in (eRe)^{\otimes_E n}$.
\end{proof}

\fait{lemma}{radical3}{If $R$ is oriented by a set of idempotents
$e, e'_1$ and $e'_2$, then the following hold, for $i=1,2$:
\begin{enumerate}
    \item[(a)] If $eRe'_i =0$, then $A_i^{{\otimes}_E n} =
e'_iR^{{\otimes}_E n}e'_i \oplus e'_iR^{{\otimes}_E n}e \oplus
C^{{\otimes}_E n}$;
    \item[(b)] If $e'_iRe =0$, then $A_i^{{\otimes}_E n} =
e'_iR^{{\otimes}_E n}e'_i \oplus eR^{{\otimes}_E n}e'_i \oplus
C^{{\otimes}_E n}$.
\end{enumerate}}
\begin{proof} This follows from $A_i^{{\otimes}_E n} =
(e_iRe_i)^{\otimes_E n}=e_i R^{{\otimes}_E n} e_i = (e'_i
+e)R^{{\otimes}_E
n}(e'_i+e)$ and Remark \ref{rem tenseurs} (c).%
\end{proof}

%
%
\section{Hochschild (co)homology}
    \label{Hochschild (co)homology}

\subsection{Definitions and notations}%
    \label{Definitions and notations}
Let, as before,  $k$ be a field and $R$ be a $k$-algebra. For each
$n\geq 2$, consider the morphism of $R$-bimodules
\mor{b'}{R^{\otimes n+1}}{R^{\otimes n}} given by
\begin{eqnarray} \label{b'} \qquad b'(r_0 \otimes \cdots \otimes
r_{n}) & := & \sum_{i=0}^{n-1}(-1)^i(r_0\otimes \cdots\otimes
r_ir_{i+1} \otimes \cdots\otimes r_n)
\end{eqnarray}
Since $b'\circ b'=0$, we get a complex $(R^{\otimes \bullet},b')$.
The maps  \mor{h}{R^{\otimes n}}{R^{\otimes n+1}} given by
$h(r_0\otimes \cdots \otimes r_{n-1})=1\otimes r_0\otimes \cdots
\otimes r_{n-1}$ define a contracting homotopy, so $(R^{\otimes
\bullet}, b')$ is a resolution of the $R^e$-module $R$, the
so-called \textbf{Hochschild resolution}.

Let $M$ be a $R$-bimodule. If we apply the functor $M\otimes_{R^e}-$
to the Hochschild resolution of $R$, we get a complex whose $n$th
homology group is, by definition, the $n$th \textbf{Hochschild
homology group} of $R$ with coefficients in $M$, denoted ${\rm
H}_n(R,M)$. Similarly, the $n$th \textbf{Hochschild cohomology
group} ${\rm H}^n(R,M)$ is the $n$th homology group of the complex
obtained upon applying the functor $\mbox{Hom}_{R^e}(-,M)$ to the
Hochschild resolution of $R$. We simply write ${\rm H}_n(R)$ and
${\rm H}^n(R)$ in case $M=R$.

Since $k$ is a field, $R^{\otimes n-2}$ is $k$-projective and so the
module $R^{\otimes n}$ is $R^e$-projective for $n\geq 2$.  Then the
Hochschild resolution $(R^{\otimes\bullet},b')$ is a projective
resolution of $R$ over $R^e$. So in our context one may define the
Hochschild (co)homology groups of $R$ with coefficients in $M$ in
the following way:
$${\rm H}_n(R,M)=\mbox{Tor}_n^{R^e}(M, R) \mbox{ \qquad and \qquad } {\rm H}^n(R,M)=\mbox{Ext}^n_{R^e}(R, M)$$
Thus, for each $n\geq 0$ there exists an isomorphism $\phi_n : D{\rm
H} ^n(R,DM)\simeq {\rm H}_n(R,M)$ (see \cite[p.~181]{CE56}), where
$D = {\rm Hom}_k(-,k)$ is the usual duality. Moreover, the
Hochschild (co)homology groups do not depend on the projective
resolution we consider to compute them. For our purpose, it will be
necessary to deal with a different projective resolution of $_RR_R$.
In particular we will use the following lemma, from
\cite{Redondo01}.

\fait{lemma}{lem proj res}{Let $R$ be an algebra and $E$ be a
semisimple subalgebra of $R$.  The complex
\begin{eqnarray} \label{Hochschild resolution E}
\qquad (R^{\otimes_E \bullet}, b') : & \SelectTips{eu}{10}\xymatrix@C=15pt{\cdots \ar[r] &
R^{\otimes_E n+1} \ar[r]^-{b'} & R^{\otimes_E n} \ar[r]^-{b'} &
\cdots \ar[r]^-{b'} & R^{\otimes_E 2} \ar[r]^-{b'} & R \ar[r] & 0}
\end{eqnarray}
where $b'$ is given by the formula (\ref{b'}) in which we replace
the tensor products over $k$ by tensor products over $E$, is a
projective resolution of $R$ as $R^e$-module.}

Hence the groups ${\rm H}_n(R,M)$ (or ${\rm H}^n(R,M)$) are the
homology groups of the complex obtained when we apply the functor
$M\otimes_{R^e}-$ (or $\mbox{Hom}_{R^e}(-,M)$, respectively) to the
projective resolution (\ref{Hochschild resolution E}).

Now, for $n\geq 2$, we have canonical isomorphisms
\begin{eqnarray}\label{tenseurs homologie}
M\otimes_{R^e}R\otimes_E R^{\otimes_E n-2}\otimes_E R& \simeq& M\otimes_{E^e} R^{\otimes_E n-2}\\
m\otimes r\otimes x\otimes r' &\mapsto &  r'mr\otimes x \nonumber
\end{eqnarray}
and
\begin{eqnarray}\label{hom cohomologie}
\mbox{Hom}_{R^e}(R\otimes_E R^{\otimes_E n-2}\otimes_E R,
M)& \simeq& \mbox{Hom}_{E^e}(R^{\otimes_E n-2}, M)\\
f &\mapsto &\overline{f}: m \mapsto f(1_R\otimes
m\otimes 1_R) \nonumber
\end{eqnarray}

Using these identifications we obtain the following (compare with \cite[(1.2)]{C89}):

\fait{proposition}{propCibils}{Let $R$ be an algebra and $E$ be a
semisimple subalgebra of $R$. Let $M$ be a $R$-bimodule.
\begin{enumerate}
\item[(a)] The Hochschild homology groups ${\rm H}_n(R,M)$ are the
homology groups of the complex
$$M\otimes_{E^e} R^{\otimes_E \bullet} \ : \ \SelectTips{eu}{10}\xymatrix@C=8pt{ \cdots
\ar[r] & M\otimes_{E^e} R^{\otimes_E 3} \ar[r]^-{\delta_2} &
M\otimes_{E^e} R^{\otimes_E 2} \ar[r]^-{\delta_1} & M\otimes_{E^e} R
\ar[r]^-{\delta_0} & M\otimes_{E^e} E \ar[r]& 0}$$ where the
differential $\delta_\bullet$ is deduced from the definition of $b'$
(\ref{Hochschild resolution E}), and the formula (\ref{tenseurs
homologie}).
\item[(b)] The Hochschild cohomology groups ${\rm H}^n(R,M)$ are the homology groups of the complex
$$\Hom_{E^e}(R^{\otimes \bullet}, M) \ : \ \SelectTips{eu}{10}\xymatrix@C=7pt{0 \ar[r] &
\Hom_{E^e}(E,M) \ar[r]^-{\delta^0} & \Hom_{E^e}(R,M)
\ar[r]^-{\delta^1} & \Hom_{E^e}(R^{{\otimes_E}2},M)
\ar[r]^-{\delta^2} & \cdots}$$ where, the differential
$\delta^\bullet$ is deduced from the definition of $b'$
(\ref{Hochschild resolution E}), and the formula (\ref{hom
cohomologie}).
\end{enumerate}}

The Yoneda product in ${\rm H}^*(R)= \bigoplus_{i\geq 0} {\rm
H}^i(R)$ is induced by a cup-product defined at the cochain level in
the complex of Proposition \ref{propCibils} (b), as follows. Given
$f \in \Hom_{E^e}(R^{\otimes_E n}, R)$ and $g \in
\Hom_{E^e}(R^{\otimes_E m}, R)$ we define the element $f \cup g$ in
$\Hom_{E^e}(R^{\otimes_E n+m}, R)$ by the rule:
$$(f \cup g)(r_1 \otimes \cdots \otimes r_{n+m}) =
f(r_1 \otimes \cdots \otimes r_n) g(r_{n+1} \otimes \cdots \otimes
r_{n+m}).$$
On the other hand, given $i \in \{1, \dots, n\}$, define $f \circ_i
g \in {\rm Hom}_{E^e}(R^{\otimes_E n+m-1},R)$ by the rule:
$$ (f \circ_i g) (r_1 \otimes \cdots \otimes r_{n+m-1})=
f(r_1 \otimes \cdots \otimes g(r_i\otimes \cdots \otimes r_{i+m-1})\otimes \cdots \otimes r_{n+m-1}).$$
Then we define the composition
$$f \circ g = \left \{
\begin{array}{ll}
0,    & \textrm{ if }n = 0{ ;}  \\ \sum_{i=1}^n (-1)^{(i-1)(m-1)}f
\circ_i g, & \textrm{ if } n > 0 {.}
\\
\end{array}
\right.$$
and finally the bracket of $f$ and $g$ as $[f,g] = f \circ g -
(-1)^{(n-1)(m-1)} g \circ f.$ As the cup-product, this bracket
induces an operation at the cohomology level. These  two structures
are related: in fact ${\rm H}^*(R)=\bigoplus_{i\geq 0} {\rm
H}^i(R)$, refered as the \textbf{Hochschild algebra} over $R$ in the
sequel, is a Gerstenhaber algebra (see \cite{G63}). That is ${\rm
H}^*(R)$ is a graded $k$-vector space endowed with a product $\cup$
which makes it a graded commutative algebra, and a Lie bracket $[-,
-]$ of degree $-1$ that makes of ${\rm H}^*(R)$ a graded Lie
algebra, and such that $[x, y\cup z]=[x,y]\cup
z+(-1)^{(|x|-1)|y|}y\cup [x,z],$ where $|x|$ and $|y|$ are the
degrees of $x$ and $y$, respectively.

\subsection{Long exact sequence of Hochschild homology groups}%
    \label{Long exact sequence of Hochschild homology groups}%

Let $R$ be an algebra, oriented by a set of idempotents $e, e'_1$
and $e'_2$ (see Definition \ref{defn oriented}). In this section, we
prove our first main result, that is the existence of a long exact
sequence relating the Hochschild homology groups of $R$, $C$, $A_1$
and $A_2$. This will be a consequence of the following proposition,
in which we freely use the identifications given in Remark \ref{rem
tenseurs} and Lemma \ref{radical2}.

\fait{proposition}{prop suite exacte 1}{Let $R$ be an algebra,
oriented by three idempotents $e, e'_1$ and $e'_2$. Moreover, let
$C=eRe$, $A_1=e_1Re_1$ and $A_2=e_2Re_2$ be as in Definition
\ref{defn oriented}.  Let $M$ be an $R$-bimodule. For any $n\geq 0$,
the sequence
$$\SelectTips{eu}{10}\xymatrix@C=20pt{0\ar[r]&  eMe\bigotimes_{E_C^e}C^{\otimes_{E_C} n} \ar[r]^-{\alpha^n} &\displaystyle  \bigoplus_{i=1}^2
e_iMe_i\otimes_{E_i^e}\!A_i^{\otimes_{E_i} n} \ar[r]^-{\beta^n} &
M \bigotimes_{E^e}R^{\otimes_E n} \ar[r] & 0} $$
is exact in \emph{mod}$E^e$, where $\alpha^n=(\alpha_1^{n},
\alpha_2^{n})^t$, with $\alpha_i^{n}(eme \otimes exe)=
e_iemee_i\otimes e_iexee_i$, and $\beta^{n}=(\beta_1^{n},
\beta_2^{n})$, with $\beta_i^{n}(e_ime_i\otimes e_ixe_i)= (-1)^i
(e_ime_i\otimes e_ixe_i)$, for $i=1,2$.}
\begin{proof}
Since $R$ is oriented, we have three cases
to consider, these are (\ref{case 1})-(\ref{case 3}). We only prove
the exactness for the case (\ref{case 1}) since the proofs for the
other cases are similar. So, assume that $e'_1Re'_2 = e'_2Re'_1 = 0$
and also $eRe'_1=eRe'_2=0$. Using Remark \ref{rem tenseurs} (c), we
have $e'_1Xe'_2=e'_2Xe'_1=0$, and similarly $e'_1Xe=e'_2Xe=0$ for
any $R$-bimodule $X$.

Using $1_R=e'_1 + e + e'_2$, we have the decomposition :
\begin{eqnarray*}
R^{\otimes_E n} = e'_1R^{\otimes_E n}e'_1 \oplus eR^{\otimes_E
n}e'_1 \oplus eR^{\otimes_E n}e \oplus eR^{\otimes_E n}e'_2 \oplus
e'_2R^{\otimes_E n}e'_2
\end{eqnarray*}
and so
\begin{eqnarray*}
M\otimes_{E^e}R^{\otimes_E n} & = &
(M\otimes_{E^e}e'_1R^{\otimes_E n}e'_1) \oplus
(M\otimes_{E^e}eR^{\otimes_E n}e'_1) \\ && \oplus
(M\otimes_{E^e}eR^{\otimes_E n}e)
\\ && \oplus
(M\otimes_{E^e}eR^{\otimes_E n}e'_2) \oplus (M\otimes_{E^e}e'_2R^{\otimes_E n}e'_2)\\ %
&\simeq & (e_1Me_1\otimes_{E_1^e}e'_1R^{\otimes_E n}e'_1) \oplus
(e_1Me_1\otimes_{E_1^e}eR^{\otimes_E n}e'_1) \\ && \oplus
(eMe\otimes_{E_C^e}C^{\otimes_{E_C} n})
\\ && \oplus (e_2Me_2\otimes_{E_2^e}eR^{\otimes_E n}e'_2) \oplus (e_2Me_2\otimes_{E_2^e}e'_2R^{\otimes_E
n}e'_2)
\end{eqnarray*}
where the isomorphism follows from immediate identification.

Similarly, for $i=1,2$, we get from Lemma \ref{radical3}, that
\begin{eqnarray*}
e_iMe_i\otimes_{E_i^e}A_i^{\otimes_{E_i} n} &\simeq&
(e_iMe_i\otimes_{E_i^e}e'_iR^{\otimes_E n}e'_i) \oplus
(e_iMe_i\otimes_{E_i^e}eR^{\otimes_E n}e'_i) \\ && \oplus
(eMe\otimes_{E_C^e}C^{\otimes_{E_C} n})
\end{eqnarray*}
With these isomorphisms and identifications in mind, we obtain the
desired short exact sequence of complexes.
\end{proof}

\fait{remark}{rem k}{\emph{%
In the above proposition, the fact that the tensor products are
taken over $E$ is crucial.  We have used Lemmata \ref{radical2} and
\ref{radical3}, which do not hold true if  tensor products are taken
over $k$.}}

We can know prove our first main result.

\fait{theorem}{thm Hochschild homology}{Let $R$ be an algebra,
oriented by three idempotents $e, e'_1$ and $e'_2$. Moreover, let
$C=eRe$, $A_1=e_1Re_1$ and $A_1=e_2Re_2$ be as in Definition
\ref{defn oriented}. For any $R$-bimodule $M$, there exists a long
exact sequence of Hochschild homology groups :
$$\SelectTips{eu}{10}\xymatrix@C=12pt@R=3pt{ \cdots \ar[r] & {\rm H}_1(C, eMe) \ar[r] &
{\rm H}_1(A_1, e_1Me_1) \oplus {\rm H}_1(A_2, e_2Me_2) \ar[r]&
{\rm H}_1(R, M) \ar[r] & \ \ \ \\ %
\ \ar[r] & {\rm H}_0(C, eMe)\ar[r] & {\rm H}_0(A_1, e_1Me_1) \oplus {\rm H}_0(A_2,
e_2Me_2) \ar[r] & {\rm H}_0(R, M) \ar[r] & 0}$$}
\begin{proof}
The result follows from the above
Proposition and the observation that ${\alpha}^{n}$ and
${\beta}^{n}$ commute with the differential $\delta_n$ of
Proposition \ref{propCibils} (a), giving rise to a short exact
sequence of complexes.
\begin{eqnarray*}
\SelectTips{eu}{10}\xymatrix@C=15pt{0\ar[r]& eMe\otimes_{E_C^e}
C^{\otimes_{E_C} \bullet} \ar[r]^-{\alpha^\bullet} & \displaystyle
\bigoplus_{i=1}^2 e_iMe_i\otimes_{E_i^e} A_i^{\otimes_{E_i} \bullet}
\ar[r]^-{\beta^\bullet}  & M\otimes_{E^e} R^{\otimes_E \bullet}
\ar[r] & 0}
\end{eqnarray*}
whose homology groups are the desired Hochschild homology groups.
\end{proof}

\subsection{Long exact sequence of Hochschild cohomology groups}%
    \label{Long exact sequence of Hochschild cohomology groups}%

In this section, we show an analogue to Theorem %
\ref{thm Hochschild homology} for Hochschild cohomology groups. In
addition, we show that the maps involved in the long exact sequence
are compatible with the Yoneda product and the Lie bracket of the
involved Hochschild algebras.

At this point, we stress that one can deduce a long exact sequence
of Hochschild cohomology groups by applying the isomorphism $\phi_n
: D{\rm H}^n(R, DM)\simeq {\rm H}_n(R,M)$
on the long exact sequence of Theorem %
\ref{thm Hochschild homology}. We will however avoid this approach
since we need a clear description of the morphisms to verify their
compatibility with the Yoneda products and Lie brackets.  Their
knowledge will also be important for Theorem \ref{cyclic connes
coh}. Nevertheless, since our approach is similar to the one we
have used for Hochschild homology, we skip some details, for instance the proof of the
following proposition, which is an analogue to Proposition
\ref{prop suite exacte 1}.

\fait{proposition}{prop suite exacte 2}{Let $R$ be an algebra,
oriented by three idempotents $e, e'_1$ and $e'_2$. Moreover, let
$C=eRe$, $A_1=e_1Re_1$ and $A_1=e_2Re_2$ be as in Definition
\ref{defn oriented}.  Let $M$ be an $R$-bimodule. For any $n\geq 0$,
the sequence
\begin{eqnarray*}
\SelectTips{eu}{10}\xymatrix@C=7pt{0\ar[r]& \Hom_{E^e}(R^{\otimes_{E}n},M)
\ar[r]^-{\alpha^n} &\displaystyle \bigoplus_{i=1}^2
\Hom_{E_i^e}(A_i^{\otimes_{E_i}n},e_iMe_i) \ar[r]^-{\beta^n} &
\Hom_{E_C^e}(C^{\otimes_{E_C}n},eMe) \ar[r] & 0}
\end{eqnarray*}
is exact in \emph{mod}$E^e$, where $\alpha^n=(\alpha_1^{n},
\alpha_2^{n})^t$, with $\alpha_i^{n}(f)(r_1\otimes\cdots\otimes
r_n)=f(e_ir_1\otimes\cdots\otimes r_ne_i)$ and
$\beta^{n}=(\beta_1^{n}, \beta_2^{n})$, with
$\beta_i^{n}(f)(r_1\otimes\cdots\otimes r_n)=(-1)^i
f(er_1\otimes\cdots\otimes r_ne)$, for $i=1,2$.}

We can now prove our second main result.

\fait{theorem}{thm Hochschild cohomology}{Let $R$ be an algebra,
oriented by three idempotents $e, e'_1$ and $e'_2$. Moreover, let
$C=eRe$, $A_1=e_1Re_1$ and $A_1=e_2Re_2$ be as in Definition
\ref{defn oriented}. For any $R$-bimodule $M$, there exists a long
exact sequence of Hochschild cohomology groups:
$$\SelectTips{eu}{10}\xymatrix@C=12pt@R=3pt{ 0 \ar[r] & {\rm H}^0(R, M)
\ar[r]^-{\alpha^0} & {\rm H}^0(A_1, e_1Me_1) \oplus {\rm H}^0(A_2,
e_2Me_2) \ar[r]^-{\beta^0} &
{\rm H}^0(C, eMe) \ar[r] & \ \ \ \\ %
\ \ar[r] & {\rm H}^1(R, M) \ar[r]^-{\alpha^1} & {\rm H}^1(A_1,
e_1Me_1) \oplus {\rm H}^1(A_2, e_2Me_2) \ar[r]^-{\beta^1} &
{\rm H}^1(C, eMe) \ar[r] & \cdots}$$
In addition, if $M=R$, then the induced maps
\mor{\alpha_i^{n}}{{\rm H}^n(R)}{{\rm H}^n(A_i)} and
\mor{\beta_i^{n}}{{\rm H}^n(A_i)}{{\rm H}^n(C)} are compatible
with the Yoneda products and the Lie brackets.}
\begin{proof}
It is easily verified that the maps
${\alpha}^{n}$ and ${\beta}^{n}$ of the above Proposition commute
with the differential $\delta^n$ of Proposition \ref{propCibils}
(b), giving rise to a short exact sequence of complexes
\begin{eqnarray*}
\SelectTips{eu}{10}\xymatrix@C=15pt{0\ar[r]& \Hom_{E^e}(R^{\otimes \bullet}, M)
\ar[r]^-{\alpha^\bullet} &\displaystyle  \bigoplus_{i=1}^2 \Hom_{E_i^e}(A_i^{\otimes \bullet},
e_iMe_i) \ar[r]^-{\beta^\bullet} & \Hom_{E_C^e}(C^{\otimes \bullet}, eMe) \ar[r]
& 0}
\end{eqnarray*}
whose homology groups are the desired Hochschild cohomology groups.
This shows the existence of the long exact sequence.  In order to
prove the second part of the statement, we need to show that
$\alpha^\bullet=(\alpha_1^\bullet, \alpha_2^\bullet)^t$ is
compatible with the cup-product and the Lie bracket if $M=R$. For
the cup-product, this follows from the fact that if $x \in
A_i^{\otimes_{E_i} (n+m)}$, then $x=e_ixe_i$ and
$$(\alpha_i^n(f)) \cup (\alpha_i^m(g)) (x)= (f \cup g) (x) = \alpha_i^{n+m}(f \cup g) (x),$$
for $f\in \Hom_{E^e}(R^{\otimes_E n}, R)$ and $g\in
\Hom_{E^e}(R^{\otimes_E m}, R)$. A similar argument applies to
$\beta_i^\bullet$.
For the Lie bracket, let $f\in \Hom_{E^e}(R^{\otimes_E n}, R)$,
$g\in \Hom_{E^e}(R^{\otimes_E m}, M)$ and
$x_1\otimes \cdots \otimes x_{n+m-1} \in A_i^{\otimes n+m-1}$. Then \\ \ \\ %
$\alpha_i^{n+m-1} (f \circ_j g) (x_1 \otimes \cdots \otimes x_{n+m-1})$

{\setlength\arraycolsep{2pt}
\begin{eqnarray*}
\ & = & (f \circ_j g) (e_ix_1 \otimes \cdots \otimes x_{n+m-1}e_i)
\\ \ & = & f(e_ix_1 \otimes \cdots \otimes g(x_j\otimes \cdots
\otimes x_{j+m-1})\otimes \cdots \otimes x_{n+m-1}e_i) \\ \ & = &
f(e_ix_1 \otimes \cdots \otimes g(e_ix_j\otimes \cdots \otimes
x_{j+m-1}e_i)\otimes \cdots \otimes x_{n+m-1}e_i) \\\ & = &
(\alpha_i^n (f)) \circ_j (\alpha_i^m (g)) (x_1 \otimes \cdots \otimes
x_{n+m-1}).
\end{eqnarray*}}
The result follows by linearity. Similarly, $\beta_i^\bullet$ is
compatible with the Lie bracket.
\end{proof}

\subsection{The case where $C$ is a core algebra}
    \label{The case where C is a core algebra}

This section is devoted to the study of a particular class of
oriented bound quiver algebras.

More precisely, let $R=kQ_R/I_R$ where $Q_R$ is a finite quiver and
$I_R$ is an admissible ideal in $kQ_R$.  Moreover, the idempotents
$e$, $e'_1$ and $e'_2$ are sums of primitive idempotents
$\varepsilon_x$ corresponding to some vertices $x$ in $Q_R$, and
$Q_C$, $Q_{A_1}$ and $Q_{A_2}$ are the full and convex subquivers of
$Q_R$ generated by these vertices.  Finally, we have $I_C=kQ_C\cap
I_R$, $I_1=kQ_{A_1}\cap I_R$ and $I_2=kQ_{A_2}\cap I_R$. Thus,
$C=kQ_C/I_C$, $A_1=kQ_{A_1}/I_1$ and $A_2=kQ_{A_2}/I_2$ are also
bound quiver algebras. For instance, the oriented algebra of
Examples \ref{examples} (a) below is like so.

We show that if $R$ is an oriented algebra as above and $C$
satisfies some specific conditions, then the long exact sequence
of Theorem \ref{thm Hochschild cohomology} results in concrete
formul\ae \ to compute the Hochschild cohomology groups of $R$.
The conditions we impose on $C=kQ_C/I_C$ are the following:%
\begin{enumerate}
\item[(1)] $Q_C$ contains no oriented
cycles, and
\item[(2)] ${\rm H}^n(C)=0$ for $n\geq 1$.
\end{enumerate}
In the sequel, an algebra satisfying the above two conditions is
called a \textbf{core algebra}. There exists many important classes
of algebras which are core algebras, mostly related to trees.

\fait{examples}{core algebras}{\emph{%
\begin{enumerate}
  \item [(a)] Let $A=kQ/I$, where $I$ is a  two-sided \textbf{monomial} ideal and $Q$
is connected (for instance $A=kQ$, a hereditary connected algebra). Then
$A$ is core if and only if $Q$ is a tree, by \cite[(1.6),(2.2)]{Hap89}
\item[(b)] We recall from \cite{Hap89} that an algebra $A$ is called \textbf{piecewise
hereditary of type $H$} if there exists a hereditary algebra $H$
such that the bounded derived categories of finitely generated
$A$-modules and $H$-modules are equivalent as triangulated
categories.
In this case, it follows from \cite[(4.2)]{Hap89} that ${\rm H}^n(A)$ is
isomorphic to ${\rm H}^n(H)$ for $n\geq 0$.  Moreover, since, by
\cite[(IV.1.10)]{H88} the quiver of any piecewise hereditary
algebra has no oriented cycles, then it follows from Example (a)
that $A$ is a core algebra if and only if the quiver of $H$ is a
tree.
\item[(c)] Following \cite{GPPRT99}, for an algebra $A=kQ/I$,
we say that ${\rm H}^n(A)$ \textbf{strongly vanishes} (see also \cite[2.2, 4.1]{Skow92} if ${\rm
H}^n(B)=0$ for every full convex subcategory $B$ of $A$. Also, we
say that $A$ has \textbf{convex projectives} (or \textbf{convex
injectives}) if for every vertex $x$ of $Q$, the set of vertices
$\{y \ | \ \varepsilon_xA\varepsilon_y\neq 0\}$ (or $\{y \ | \
\varepsilon_yA\varepsilon_x\neq 0\}$, respectively) is convex in
$Q$. With these definitions, the following holds true by
\cite{GPPRT99} : $A$ is a core algebra provided %
\begin{enumerate} \item[(i)] $A$ is schurian and ${\rm H}^1(A)$ strongly
vanishes, or \item[(ii)] $Q$ has no oriented cycles and $A$ is a
tame algebra with convex projectives and convex injectives such that
${\rm H}^1(A)$ and ${\rm H}^2(A)$ strongly vanish. \end{enumerate}
\end{enumerate}}}

We continue with the following lemma concerning the center of an
algebra, whose proof is left as an exercise. Recall that the $0$th Hochschild cohomology group ${\rm
H}^0(A)$ of an algebra $A$ is isomorphic to its center $Z(A)$.

\fait{lemma}{center}{Let $A=kQ/I$ be a connected algebra,
$\varepsilon_i$ be the idempotent corresponding to the vertex $i$,
and $\omega$ be an element of the center $Z(A)$ of $A$. Then,
$\omega = \sum_{i\in Q_0}\omega_{i} + \lambda \cdot 1_A$, where
$\lambda \in k$ and $\omega_{i}\in \varepsilon_i \rad A
\varepsilon_i$ for each $i$.}
%

\fait{remark}{component}{\emph{It follows from the previous lemma that if
$C=kQ_C/I_C$ is an algebra such that $Q_C$ contains no oriented
cycles, then dim$_k {\rm H}^0(C)$ equals the number of connected
components of $Q_C$.}}

We can now prove the main result of this section.

\fait{proposition}{Hochart}{Let $R=kQ_R/I_R$ be a bound quiver
algebra, oriented by three idempotents $e, e'_1$ and $e'_2$.
Moreover, let $C=kQ_{C}/I_C$, $A_1=kQ_{A_1}/I_1$ and
$A_2=kQ_{A_2}/I_2$. If $C=kQ_C/I_C$ is a core algebra, then:
\begin{enumerate}
\item [(a)] $\emph{dim}_k {\rm H}^0(R) = \emph{dim}_k {\rm H}^0(A_1) +
\emph{dim}_k {\rm H}^0(A_2) - 1$; %
\item[(b)] $\emph{dim}_k{\rm H}^1(R) = \emph{dim}_k{\rm H}^1(A_1) +
\emph{dim}_k{\rm H}^1(A_2) + \emph{dim}_k{\rm H}^0(C)-1$;%
\item[(c)] ${\rm H}^i(R) \simeq {\rm H}^i(A_1) \oplus {\rm H}^i(A_2)$ for $i
\geq 2$.
\end{enumerate}}
\begin{proof}
(a). It is sufficient to show that $\mbox{dim}_k Z(R) = \mbox{dim}_k
Z(A_1) + \mbox{dim}_k Z(A_2)-1$. For $j=1,2$, let $\mathbb{B}_j$ be
a $k$-basis of $Z(A_j)$ containing $e_j$, the unit element of $A_j$.
We claim that $\mathbb{B}= (1_R \cup \mathbb{B}_1 \cup \mathbb{B}_2)
\setminus (\{e_1\} \cup \{e_2\})$ is a $k$-basis of $Z(R)$. Since
this set is linearly independent by Lemma \ref{center}, it remains
to show that it generates $Z(R)$. Let $\omega$ be an element of
$Z(R)$. By Lemma \ref{center}, we have $\omega = \sum_{i\in
{(Q_R)}_0}\omega_{i} + \lambda \cdot 1_R$, where $\lambda\in k$ and
$\omega_{i}\in \varepsilon_i\mbox{rad}R\varepsilon_i$ for each $i\in
{(Q_R)}_0$. Moreover, since $Q_C$ has no oriented cycles, we have
$\omega = \sum_{i\in {(Q_R)}_0\setminus {(Q_C)}_0}\omega_{i} +
\lambda \cdot 1_R$. Therefore, we can write without ambiguity
$\omega = \rho_1 + \rho_2 + \lambda \cdot 1_R$, where $\rho_j
=\sum_{i\in {(Q_{A'_j})}_0}\omega_i$, $j=1,2$. It then remains to
show that $\rho_j\in Z(A_j)$ for $j=1,2$ \ since, by construction,
$\rho_j\in Z(A_j)$ if and only if $\rho_j$ belongs to the $k$-vector
space generated by $\mathbb{B}_j\setminus\{e_j\}$. Let $\eta_1$ be a
path in $Q_{A_1}$.  We have :
$$ \begin{array}{rl}
  \rho_1\eta_1 + \lambda\eta_1 & = \rho_1\eta_1 + 0 + \lambda\eta_1 = \rho_1\eta_1 + \rho_2\eta_1 + \lambda \cdot 1_R \eta_1
   = \omega \eta_1 \\
   & = \eta_1\omega = \eta_1\rho_1 + \eta_1\rho_2 + \lambda\eta_1 \cdot 1_R = \eta_1\rho_1 + 0 + \lambda\eta_1\\
   & = \eta_1\rho_1 + \lambda\eta_1
\end{array}$$
whence $\rho_1\eta_1 = \eta_1\rho_1$.  Therefore, $\rho_1\in
Z(A_1)$. Since, similarly, $\rho_2\in Z(A_2)$, $\omega$ belongs to
the $k$-vector space generated by $\mathbb{B}$ and $\mathbb{B}$ is
a $k$-basis of $Z(R)$. In particular, $\mbox{dim}_k Z(R) =
\mbox{dim}_k Z(A_1) + \mbox{dim}_k Z(A_2)-1$, and this proves (a).

(b). Since ${\rm H}^i(C) = 0$ for each $i \geq 1$, it follows from
Theorem \ref{thm Hochschild cohomology} that we have the following
exact sequence :
$$\SelectTips{eu}{10}\xymatrix@C=12pt@R=12pt{ 0 \ar[r] & {\rm H}^0(R) \ar[r] &
{\rm H}^0(A_1) \oplus {\rm H}^0(A_2) \ar[r] & {\rm H}^0(C) \ar[r] &
{\rm H}^1(R)\ar[r] & {\rm H}^1(A_1) \oplus {\rm H}^1(A_2) \ar[r]& 0}$$
whence, %
$$ \begin{array}{rl}
  \mbox{dim}_k {\rm H}^1(R) = & \mbox{dim}_k {\rm H}^1(A_1) + \mbox{dim}_k {\rm H}^1(A_2) + \ \mbox{dim}_k {\rm H}^0(C) \\
  &  + \mbox{dim}_k {\rm H}^0(R) -
\mbox{dim}_k {\rm H}^0(A_1) - \mbox{dim}_k {\rm H}^0(A_2)
\end{array}$$
But $\mbox{dim}_k {\rm H}^0(R) - \mbox{dim}_k {\rm H}^0(A_1) -
\mbox{dim}_k {\rm H}^0(A_2)= -1$ by (a).

(c). This follows directly from Theorem \ref{thm Hochschild
cohomology}.
\end{proof}

\subsection{Examples and comparison}
    \label{Examples and comparison}

We now compute some examples and compare the sequence of Theorem
\ref{thm Hochschild cohomology} with other ones existing in the
literature. First, let
$$R= \left(
\begin{array}{cc} A & 0  \\ _BM_A & B
\end{array} \right)$$
where $A$ and $B$ are arbitrary finitely generated $k$-algebras and
$M$ is a finitely generated $B$-$A$-bimodule. Then $R$ is a
$k$-algebra under the usual matrix operations.

In the eighties, D.~Happel proved in \cite{Hap89} the existence of
a long exact sequence relating the Hochschild cohomology groups of
$R$, $A$ and $B$ in the particular case where $B=k$. Over the
years, several variations and generalizations of this sequence
have been obtained.  The most recent long exact sequences now
apply for an arbitrary algebra $B$ (see \cite{Ci00, MP00, GS02,
GMS03, BG04}).

\fait{theorem}{Michelena-Platzeck}{Let $R$ be as above, $e_A$ and
$e_B$ are the unit elements of $A$ and $B$ and let $X$ be a
$R$-$R$-bimodule.  Then, there exists a long exact sequence
$$\SelectTips{eu}{10}\xymatrix@C=8pt@R=3pt{ 0 \ar[r] & {\rm H}^0(R, X) \ar[r] &
{\rm H}^0(A, e_AXe_A) \oplus {\rm H}^0(B, e_BXe_B) \ar[r] &
\Ext_{B\otimes_kA^{op}}^0(M, e_BXe_A) \ar[r] & \ \ \ \\ %
\ \ar[r] & {\rm H}^{1}(R, X)\ar[r] & {\rm H}^1(A, e_AXe_A) \oplus
{\rm H}^1(B, e_BXe_B) \ar[r] & \Ext_{B\otimes_kA^{op}}^1(M, e_BXe_A)
\ar[r] & \cdots}$$}

The following examples show the difference between the exact
sequence of Theorem \ref{thm Hochschild cohomology} and the above
sequence.

\fait{examples}{examples}{ \emph{ %
\begin{enumerate} %
\item[(a)] Let $R$ be the bound quiver algebra given by the quiver%
$$\SelectTips{eu}{10}\xymatrix@R=8pt@C=15pt{5 \ar[d] \ar[dr] \ar[drr] & & 6 \ar[d] \ar[dl] \ar[dll]\\
2\ar[dr]& 3 \ar[d] &4\ar[dl]\\
&1& }$$ bounded by the ideal generated by all paths of length two.
Let $e=\varepsilon_1+\varepsilon_2+ \varepsilon_3 + \varepsilon_4$,
$e'_1=\varepsilon_5$, and $e'_2=\varepsilon_6$ be the three
idempotents orienting $R$. By \cite[(p.~96)]{Ci98},
$$\qquad {\setlength\arraycolsep{4pt}\begin{array}{cc}
  {\rm H}^i(A_1)\simeq{\rm H}^i(A_2) \simeq \left\{%
{\setlength\arraycolsep{2pt}\begin{array}{ll}
  k, & \mbox{if } i=0; \\
  k^2, & \mbox{if } i=1; \\
  0, & \mbox{otherwise};
\end{array}}\right. \mbox{ and }& %
{\rm H}^i(R) \simeq \left\{%
{\setlength\arraycolsep{2pt}\begin{array}{ll}
  k, & \mbox{if } i=0; \\
  k^4, & \mbox{if } i=1; \\
  0, & \mbox{otherwise}.
\end{array}}\right.
\end{array}}
$$%
Moreover, since the quiver of $C$ is a tree, ${\rm H}^0(C)=k$ and
${\rm H}^i(C)=0$ for $i\geq 1$.  By Theorem \ref{thm Hochschild
cohomology} and Proposition \ref{Hochart}, there exists an exact
sequence : \\
$\SelectTips{eu}{10}\xymatrix@R=10pt@C=8pt{%
^{0} \ar[r] & ^{{\rm H}^0(R)} \ar[r] & ^{{\rm H}^0(A_1)\oplus
{\rm H}^0(A_2)} \ar[r] & ^{{\rm H}^0(C)} \ar[r] & ^{{\rm H}^1(R)} \ar[r]
& ^{{\rm H}^1(A_1)\oplus {\rm H}^1(A_2)} \ar[r] & ^{{\rm H}^1(C)} \ar[r]
& ^{\cdots} \\
^{0} \ar[r] & *+[F--]{^{k}} \ar[r] \ar@{=}[u] & ^{k\oplus k} \ar[r]
\ar@{=}[u] & ^{k} \ar[r] \ar@{=}[u] & *+[F--]{^{k^4}} \ar[r]
\ar@{=}[u] & ^{k^2\oplus k^2} \ar[r] \ar@{=}[u] & ^{0} \ar[r]
\ar@{=}[u] & ^{\cdots}}$\\
On the other hand, if we respectively denote by $S_x$ and $P_x$
the simple module and the projective module associated to the
vertex $x$ in the quiver, one can see $R$ as the one-point
extension $$A_1[M]=\left(\begin{array}{cc}
  A_1 & 0 \\
  M & k
\end{array}\right),$$
where $M=\mbox{rad}P_6 = S_2\oplus S_3 \oplus S_4$ is a semisimple
module without self-extensions.  In this case, the long exact
sequence of Theorem \ref{Michelena-Platzeck} is given by :\\
$\SelectTips{eu}{10}\xymatrix@R=10pt@C=8pt{%
^{0} \ar[r] & ^{{\rm H}^0(R)} \ar[r] & ^{{\rm H}^0(A_1)\oplus k}
\ar[r] & ^{\Hom_{A_1}(M,M)} \ar[r] & ^{{\rm H}^1(R)} \ar[r] &
^{{\rm H}^1(A_1)} \ar[r] &
^{\Ext^1_{A_1}(M,M)} \ar[r] & ^{\cdots} \\ %
^{0} \ar[r] & ^{k} \ar[r] \ar@{=}[u] & ^{k\oplus k} \ar[r]
\ar@{=}[u] & ^{k^3} \ar[r] \ar@{=}[u] & ^{k^4} \ar[r] \ar@{=}[u] &
^{k^2} \ar[r] \ar@{=}[u] & ^{0} \ar[r] \ar@{=}[u] & ^{\cdots}}$
\\ %
The two long exact sequences are clearly different.
\item[(b)] Let $R$ be the algebra given by the quiver
$$\SelectTips{eu}{10}\xymatrix{1 \ar@(ul, ur)[]^\varepsilon
& 2 \ar[l]_\delta & 3 \ar[l]_\gamma & 4 \ar[l]_\beta
\ar@(ul, ur)[]^\alpha}$$ with the relations $\alpha^2 =
\alpha\beta= \beta\gamma = \delta\varepsilon =
\varepsilon^2=0$.  Let $e=\varepsilon_3$,
$e'_1=\varepsilon_1 + \varepsilon_2$ and
$e'_2=\varepsilon_4$. Then the conditions (\ref{case 0})
and (\ref{case 3}) of Definition \ref{defn oriented} are
satisfied, and so $R$ is oriented. A quick application
of \cite[(p.~96)]{Ci98} and Theorem \ref{Michelena-Platzeck}
gives:%
$$ {\rm H}^i(A_1)\simeq{\rm H}^i(A_2) \simeq \left\{%
\begin{array}{ll}
  k^2, & \mbox{ if } i=0; \\
  k, & \mbox{ otherwise }.
\end{array}
\right.$$
Moreover, ${\rm H}^0(C)=k$ and ${\rm H}^i(C)=0$ for $i\geq 1$.
By Theorem \ref{thm Hochschild cohomology} and Proposition
\ref{Hochart}, there exists an exact sequence \\
$\SelectTips{eu}{10}\xymatrix@R=10pt@C=8pt{%
^{0} \ar[r] & ^{{\rm H}^0(R)} \ar[r] & ^{{\rm H}^0(A_1)\oplus
{\rm H}^0(A_2)} \ar[r] & ^{{\rm H}^0(C)} \ar[r] & ^{{\rm H}^1(R)} \ar[r]
& ^{{\rm H}^1(A_1)\oplus {\rm H}^1(A_2)} \ar[r] & ^{{\rm H}^1(C)}
\ar[r] & ^{\cdots} \\ %
^{0} \ar[r] & *+[F--]{^{k^3}} \ar[r] \ar@{=}[u] & ^{k^2\oplus k^2}
\ar[r] \ar@{=}[u] & ^{k} \ar[r] \ar@{=}[u] & *+[F--]{ ^{k^2}} \ar[r]
\ar@{=}[u] & ^{k\oplus k} \ar[r] \ar@{=}[u] & ^{0} \ar[r] \ar@{=}[u]
& ^{\cdots}}$\\
from which we obtain that ${\rm H}^1(R) \simeq k^2$ and ${\rm H}^i(R)
\simeq {\rm H}^i(A_1)\oplus {\rm H}^i(A_2)\simeq k\oplus k$ for $i\geq
2$. %
On the other hand, one can see $R$ as the following matrix algebra
$$R = \left(\begin{array}{cc}
  A_1 & 0 \\
  M & A'_2
\end{array}\right),$$
where $A'_2$ is the full subcategory generated by the vertex 4 and
$M$ is the $A'_2$-$A_1$-bimodule given by the arrow $\beta$. In
order to apply the sequence of Theorem \ref{Michelena-Platzeck},
one needs to compute $\Ext^i_{A'_2 \otimes {A_1}^{op}}(M,M)$ for
$i\geq 0$. A direct computation leads to the fact that ${A'_2}
\otimes {A_1}^{op}$ is given by the quiver
$$\SelectTips{eu}{10}\xymatrix{4\otimes1 \ar@(ul, ur)[]^{4\otimes\varepsilon}
\ar@(dl, dr)[]_{\alpha\otimes 1} \ar[r]^{4\otimes\delta} &
4\otimes 2 \ar@(dl, dr)[]_{\alpha\otimes 2}
\ar[r]^{4\otimes\gamma} & 4\otimes 3 \ar@(dl, dr)[]_{\alpha\otimes
3}}$$ bound by the relations induced by those of $R$, that is
$(\alpha\otimes 1)^2=(\alpha\otimes 2)^2=(\alpha\otimes
3)^2=(4\otimes
\varepsilon)^2=(4\otimes\varepsilon)(4\otimes\delta)=0$,
$(\alpha\otimes 1)(4\otimes\delta) =
(4\otimes\delta)(\alpha\otimes 2)$ and $(\alpha\otimes
2)(4\otimes\gamma) = (4\otimes \gamma)(\alpha\otimes 3)$.
Moreover, the ${A'_2} \otimes {A_1}^{op}$-module $M$ is the simple
module $S_{4\otimes 3}$.  However, since $(\alpha\otimes 3)^2=0$,
it follows from \cite[(1.4)]{LM04} that $\Ext^i_{{A'_2}\otimes
{A_1}^{op}}(M,M)\neq 0$ for $i\geq 1$. Hence, the long exact
sequence of Theorem \ref{Michelena-Platzeck} cannot efficiently be
used to compute the Hochschild cohomology groups of $R$. At least,
it cannot be used to show that ${\rm H}^i(R) \simeq
{\rm H}^i(A_1)\oplus {\rm H}^i(A_2)\simeq k\oplus k$ for $i\geq 2$.
\end{enumerate}}}
%

%
%
\section{Cyclic (co)homology}
    \label{Cyclic (co)homology}


The strong relation between cyclic (co)homology and Hochschild
(co)homology theories brings up the natural question whether
oriented algebras give rise to long exact sequences of cyclic
(co)homology groups. In what follows, we show that the answer is
positive, and moreover that these long exact sequences are
compatible with Connes' long exact sequences.

The first part of this section is devoted to basic definitions and
tools concerning cyclic (co)homology, including the Connes' long
exact sequences.  Then we consider the oriented algebras again and
prove our main results.

We refer to \cite{L92} and \cite{BII98} for more details on cyclic
(co)homology theory.

\subsection{Definitions and notations}%
    \label{Definitions and notations : cyclic}
We first recall the following general constructions. A
\textbf{bicomplex} $C_{\bullet \bullet}$ is a collection of modules
$C_{pq}$, with $p,q\in\mathbb{Z}$, together with a horizontal
differential \mor{d^h}{C_{p,q}}{C_{p-1,q}} and a vertical
differential \mor{d^v}{C_{p,q}}{C_{p,q-1}}
%
%
satisfying $d^vd^v=d^hd^h=d^vd^h+d^hd^v=0$. Now, given a bicomplex
$C_{\bullet \bullet}$ in the first quadrant, that is $C_{p,q}=0$
provided $p<0$ or $q<0$,  define, for $n\in \mathbb{Z}$, the module
$(\mbox{Tot} C_{\bullet \bullet})_n=\bigoplus_{p+q=n}C_{p,q}$, and,
further, the complex
$$\SelectTips{eu}{10}\xymatrix{\cdots \ar[r] & (\mbox{Tot} C_{\bullet \bullet})_{n+1} \ar[r]^d &
(\mbox{Tot} C_{\bullet \bullet})_n \ar[r]^d & (\mbox{Tot} C_{\bullet
\bullet})_{n-1} \ar[r] & \cdots}$$
where $d=d^h+d^v$.  This complex is called the \textbf{total
complex} of the bicomplex $C_{\bullet \bullet}$. Then it is a
routine verification that a short exact sequence of bicomplexes in
the first quadrant yields a short exact sequence of their associated
total bicomplexes.

Now, let $R$ be an algebra. For each $n\geq 2$, consider the
morphism of $R$-bimodules \mor{b}{R^{\otimes n}}{R^{\otimes n-1}}
given by
\begin{eqnarray}\label{b}
\qquad b(r_0\otimes \cdots \otimes r_{n-1}) & := &
\sum_{i=0}^{n-2}(-1)^i(r_0\otimes r_1\otimes \cdots\otimes
r_ir_{i+1} \otimes \cdots\otimes r_{n-1}) \\
&& + (-1)^{n-1}(r_{n-1}r_0\otimes r_1 \otimes \cdots \otimes
r_{n-2}) \nonumber
\end{eqnarray}
It is easily verified that $b\circ b=0$, so we get a complex $(R^{\otimes \bullet}, b)$.
%

On the other hand, for $n\geq 1$, the cyclic group
$\mathbb{Z}_{n+1}$ acts on $R^{\otimes n+1}$ by letting its
generator $t$ acts by $t(r_0\otimes\cdots\otimes r_n):= (-1)^n
(r_n\otimes r_0 \cdots\otimes r_{n-1})$. For $n=0$, let $t$ be the
identity map. So $t$ defines an endomorphism on $R^{\otimes n+1}$.
Let $N=1+t+\cdots +t^n$. Then the following diagram
\begin{eqnarray} \label{CC}
\SelectTips{eu}{10}\xymatrix{\vdots\ar[d] & \vdots\ar[d] & \vdots\ar[d] & \vdots\ar[d] \\ %
R^{\otimes 3} \ar[d]_b & R^{\otimes 3}\ar[d]_{-b'} \ar[l]_{1-t} &
R^{\otimes 3} \ar[d]_{b}\ar[l]_{N} & R^{\otimes 3}
\ar[l]_{1-t}\ar[d]_{-b'} & \ar[l]_{N} \cdots\\
R^{\otimes 2} \ar[d]_b & R^{\otimes 2}\ar[d]_{-b'} \ar[l]_{1-t} &
R^{\otimes 2} \ar[d]_{b}\ar[l]_{N} & R^{\otimes 2}
\ar[l]_{1-t}\ar[d]_{-b'} & \ar[l]_{N} \cdots\\
R & R \ar[l]_{1-t} & R\ar[l]_{N} & R \ar[l]_{1-t} & \ar[l]_{N}
\cdots}
\end{eqnarray}
where $b$ and $b'$ are as in (\ref{b}) and (\ref{b'}), is a
bicomplex in the first quadrant. This diagram is called the
\textbf{cyclic bicomplex} associated with $R$ and is denoted
$CC(R)$.

By definition, the $n$th \textbf{cyclic homology group} of $R$ is given by
$${\rm HC}_n(R):={\rm H}_n(\mbox{Tot}CC(R)),$$ that is the $n$th homology group of the total complex associated
with $CC(R)$.

In order to define the cyclic cohomology groups, we consider the
standard duality $D=\mbox{Hom}_k(-,k)$. Applying $D$ on the
bicomplex $CC(R)$ above yields the dual bicomplex $D(CC(R))$.
%
By definition, the $n$th \textbf{cyclic cohomology group} ${\rm
HC}^n(R)$ of $R$ is the $n$th homology group of the cochain complex
Tot$D(CC(R))$: $$ {\rm HC}^n(R):= {\rm H}^n(\mbox{Tot}D(CC(R))).$$

\subsection{Connes' exact sequences}%
    \label{Connes' exact sequences}
Let $CC(R)^{\{2\}}$ be the bicomplex consisting of the first two
columns of $CC(R)$.  By Killing's lemma (see \cite[p.~55]{L92} for
instance) the homology groups of Tot$CC(R)^{\{2\}}$ are isomorphic
to the homology groups of the first column, which are the Hochschild
homology group ${\rm H}_n(R)$ because $1_R\otimes b'$ becomes $b$
under the isomorphism $R\otimes_{R^e} R^{\otimes n+2} \simeq
R\otimes R^{\otimes n}$ given by \maps{r_0\otimes r_1\otimes \cdots
\otimes r_{n+2}}{r_{n+2}r_0r_1\otimes r_2\otimes \cdots \otimes
r_{n+1}}. Now, let  $(CC(R)[2,0])_{pq}=CC(R)_{p-2, q}$. Thus, the
short exact sequence of bicomplexes
$$\SelectTips{eu}{10}\xymatrix{0\ar[r] & CC(R)^{\{2\}} \ar[r] & CC(R) \ar[r] &
CC(R)[2,0] \ar[r] & 0}$$
yields a short exact sequence of total complexes, and so a long exact sequence
\begin{eqnarray*}
\qquad \SelectTips{eu}{10}\xymatrix@C=15pt{\cdots \ar[r] & {\rm H}_n(R)\ar[r] & {\rm HC}_n(R)
\ar[r]& {\rm HC}_{n-2}(R) \ar[r] & {\rm H}_{n-1}(R) \ar[r] & \cdots}
\end{eqnarray*}
relating the Hochschild homology groups and cyclic homology groups
of $R$, called the \textbf{homological Connes' exact sequence}.

Dually, let $D(CC(R))^{\{2\}}$ denote the bicomplex consisting of
the first two columns of $D(CC(R))$. Since $D$ is an exact functor,
Killing's lemma applies and the homology groups of
Tot$D(CC(R))^{\{2\}}$ are isomorphic to the homology groups of the
first column, which are the Hochschild cohomology groups ${\rm
H}^n(R,D(R))$ because Hom$_{R^e}(b', D(R))$ becomes $D(b)$ under the
isomorphism \mor{\phi}{\mbox{Hom}_{R^e}(R^{\otimes
n+1},D(R))}{\mbox{Hom}_k(R^{\otimes n}, k)} sending any $f\in
\mbox{Hom}_{R^e}(R^{\otimes n+1}, D(R))$ to the map $\phi(f)$
defined by $\phi(f)(r_0\otimes\cdots\otimes r_{n-1})=f(1\otimes
r_1\otimes \cdots\otimes r_{n-1}\otimes 1)(r_0)$ for
$r_0\otimes\cdots\otimes r_{n-1}\in R^{\otimes n}$. Therefore, the
short exact sequence of bicomplexes%
$$\SelectTips{eu}{10}\xymatrix{0\ar[r] & D(CC(R))[2,0] \ar[r] & D(CC(R)) \ar[r] &
D(CC(R))^{\{2\}} \ar[r] & 0}$$
(where $D(CC(R)[2,0])_{pq}=D(CC(R))_{p-2, q}$) yields a short exact
sequence of total complexes giving rise to the \textbf{cohomological
Connes' exact sequence}
\begin{eqnarray*}
\SelectTips{eu}{10}\xymatrix@C=10pt{\cdots \ar[r] &
{\rm H}^n(R,D(R))\ar[r] & {\rm HC}^{n-1}(R) \ar[r]&
{\rm HC}^{n+1}(R) \ar[r] & {\rm H}^{n+1}(R,D(R)) \ar[r] & \cdots}
\end{eqnarray*}

\subsection{Cyclic (co)homology for oriented algebras}%
    \label{Cyclic (co)homology for oriented algebras}

Let $R$ be an algebra, oriented by three idempotents $e, e'_1,
e'_2$. We also adopt the notations and terminology of Section
\ref{pullbacks}. In what follows, we give an analogues to Theorems
\ref{thm Hochschild homology} and \ref{thm Hochschild cohomology}
for the cyclic (co)homology groups.

The strategy is the same as for Hochschild (co)homology groups.
Thus, in view of Remark \ref{rem k},  it will be necessary to deal
with tensor products over $E$ rather than over $k$.  This means that
it would be nice to construct a short exact sequence relating the
bicomplexes $CC_E(C)$, $CC_E(A_1)$, $CC_E(A_2)$ and $CC_E(R)$ (where
these are obtained from the bicomplexes (\ref{CC}) by replacing the
tensor products over $k$ by tensor products over $E$), giving rise
to a long exact sequence of homology groups. However, although such
a sequence exists, this is unfortunately not sufficient because
there is no guarantee, for instance, that the total complexes of
$CC(R)$ and $CC_E(R)$ have the same homology groups. This problem
can be fixed by what is sometimes called the relative Hochschild
homology.

Let $S$ be a subring of $R$ (in our case we will take $S=E$). In
particular, $R$ is a $S$-bimodule.  Consider the $S$-subbimodule
$J^{n}$ of $R^{\otimes_S n}$ generated by all elements of the form
$(sr_1\otimes r_2\otimes \cdots \otimes r_n)-(r_1\otimes r_2\otimes
\cdots \otimes r_ns)$ for any $s\in S$ and $r_i\in R$. Given $x,
x'\in R^{\otimes_S n}$, we write
$$x\sim x' \mbox{ provided } x-x'\in J^{n}\mbox{, and define }
\widetilde{R^{\otimes_S n}}:=R^{\otimes_S n}/J^{n}. $$ %
It is straightforward to check that the differential $b$ of
(\ref{b}) is compatible with the relation $\sim$, in the sense that
if $x\sim x'$ then $b(x)\sim b(x')$.  So there is a well-defined
complex $(\widetilde{R^{\otimes_S n}}, b)$. It's $n$th homology
group is denoted $\widetilde{{\rm H}_n^S}(R)$.

Recall that an algebra $S$ is called \textbf{separable} over $k$ if
the $S$-bimodule map \mor{\mu}{S\otimes S^{op}}{S} splits. This is
equivalent to the existence of an idempotent $\varepsilon = \sum
u_i\otimes v_i \in S\otimes S^{op}$ such that $\sum u_iv_i=1$ and
$(s\otimes 1_S)\varepsilon=(1_S \otimes s)\varepsilon$ for any $s\in
S$.

Then we have the following theorem from \cite[(p.21)]{L92}.

\fait{theorem}{thm relative}{%
Let $S$ be separable over $k$. Then the canonical projection
\mor{\pi}{(R^{\otimes n}, b)}{(\widetilde{R^{\otimes_S n}}, b)}
is a quasi-isomorphism. In particular, ${\rm H}_n(R)\simeq
\widetilde{{\rm H}_n^S}(R)$ for any $n\geq 0$.}

Moreover, it is straightforward to check that if, in the diagram
(\ref{CC}), we replace the tensor products over $k$ by tensors
products over $S$, to obtain a complex $CC_S(R)$, the induced
applications $b$, $-b'$, $t$ and $N$ are compatible with the
relation $\sim$, so that we obtain a bicomplex
$\widetilde{CC_S}(R)$, refered in the sequel as the
\textbf{relative bicomplex}, given by
\begin{eqnarray} \label{CC_S}
\SelectTips{eu}{10}\xymatrix{\vdots\ar[d] & \vdots\ar[d] & \vdots\ar[d] & \vdots\ar[d] \\ %
\widetilde{R^{\otimes_S 3}} \ar[d]_b & \widetilde{R^{\otimes_S
3}}\ar[d]_{-b'} \ar[l]_{1-t} & \widetilde{R^{\otimes_S 3}}
\ar[d]_{b}\ar[l]_{N} & \widetilde{R^{\otimes_S 3}}
\ar[l]_{1-t}\ar[d]_{-b'} & \ar[l]_{N} \cdots\\
\widetilde{R^{\otimes_S 2}} \ar[d]_b & \widetilde{R^{\otimes_S
2}}\ar[d]_{-b'} \ar[l]_{1-t} & \widetilde{R^{\otimes_S 2}}
\ar[d]_{b}\ar[l]_{N} & \widetilde{R^{\otimes_S 2}}
\ar[l]_{1-t}\ar[d]_{-b'} & \ar[l]_{N} \cdots\\
\widetilde{R} & \widetilde{R} \ar[l]_{1-t} &
\widetilde{R}\ar[l]_{N} & \widetilde{R} \ar[l]_{1-t} & \ar[l]_{N}
\cdots}
\end{eqnarray}
and the map $\pi$ of Theorem \ref{thm relative} extends to a map of
bicomplexes \mor{\pi}{CC(R)}{\widetilde{CC_S}(R)}.  Moreover, with
these notations, Theorem \ref{thm relative} assures that the
restriction of $\pi$ to the columns of type $b$ is a
quasi-isomorphism.
On the other hand, since the columns of type $b'$ are exact, the
projection $\pi$ trivially induces a quasi-isomorphism on the
columns of type $b'$. Then, it follows (see \cite[Proposition
1.0.12]{L92} for instance) that the map induced by $\pi$ on the
total complexes is a quasi-isomorphism. We have just proved the
following proposition.

\fait{proposition}{prop quasi}{Let $R$ be an algebra and $S$ be a
separable subalgebra of $R$.  Then ${\rm HC}_n(R)\simeq {\rm
H}_n(\emph{Tot}\widetilde{CC_S}(R))$ for each $n\geq 0$.}

For an oriented algebra $R$, the subalgebra $S=E$ of $R$ is
separable over $k$ (take $\varepsilon=e'_1\otimes e'_1 + e\otimes
e + e'_2\otimes e'_2$), and so are the subalgebras $E_C$ of $C$,
$E_1$ of $A_1$ and $E_2$ of $A_2$.
As a consequence, to obtain a long exact sequence of cyclic
homology groups, it suffices to construct a short exact sequence
of relative bicomplexes of $R$, $C$, $A_1$ and $A_2$.  The
existence of such a sequence will follow from the next
proposition.

\fait{proposition}{prop suite exacte}{Let $R$ be an algebra,
oriented by three idempotents $e, e'_1, e'_2$.  For any $n\geq 1$,
the sequence
\begin{eqnarray*} 
\SelectTips{eu}{10}\xymatrix@C=25pt{0\ar[r]& C^{\otimes_{E_C} n}
\ar[r]^-{f^n} & \displaystyle A_1^{\otimes_{E_1} n}\bigoplus
A_2^{\otimes_{E_2} n} \ar[r]^-{g^{n}} & R^{\otimes_E n} \ar[r] &
0}
\end{eqnarray*}
is exact in \emph{mod}$E^e$, where $f^{n}=(f_1^{n}, f_2^{n})^t$,
with $f_i^{n}(er_1\otimes \cdots\otimes r_ne)= e_ier_1\otimes
\cdots\otimes r_nee_i$, and $g^{n}=(g_1^{n}, g_2^{n})$, with
$g_i^n (e_i r_1 \otimes \cdots \otimes r_n e_i) = (-1)^i (e_i r_1
\otimes \cdots \otimes r_ne_i)$, for $i=1,2$.}
\begin{proof}
Since $e_ie=e=ee_i$, the map $f^{n}$ is
trivially injective.  Also, Im$g^{n}\subseteq \mbox{Ker} f^{n}$.
Conversely, assume $x_1=\sum_{i=1}^ke_1r_{1i}\otimes \cdots \otimes
r_{ni}e_1\in A_1^{\otimes_E n}$ and
$x_2=\sum_{i=1}^le_2r'_{1i}\otimes \cdots \otimes r'_{ni}e_2\in
A_2^{\otimes_E n}$ are such that $g^{n}((x_1, x_2)) =0$.  By
definition of $g^{n}$, we get $x_1-x_2=0$, and so $x_1=x_2$.  Let
$x$ be this common value. Multiplying $x$ by $e_1$ on both sides
gives $x=x_1 = e_1x_1e_1=e_1x_2e_1 =e_1x_2e_1 =
\sum_{i=1}^{l}er'_{1i}\otimes \cdots \otimes r'_{ni}e$ since
$e_1e_2=e=e_2e_1$. So $x\in C^{\otimes_E n}$ and $f^{n}(x)=(x_1,
x_2)$, showing that Ker$g^{n}\subseteq \mbox{Im}f^{n}$. It remains
to show that $g^{n}$ is surjective.  Let $r\in R^{\otimes_E n}$.
Since $1_R=e'_1+ e + e'_2$, and $e'_1re'_2=0=e'_2re'_1$ by Remark
\ref{rem tenseurs} (c),  we have
\begin{eqnarray*}
r & = & e'_1re'_1+ e'_1re + e'_1re'_2 + ere'_1 + ere + ere'_2 +
e'_2re'_1 + e'_2re + e'_2re'_2 \\
& = & e'_1re'_1+ e'_1re + ere'_1 + ere + ere'_2 + e'_2re +
e'_2re'_2
\end{eqnarray*}
Moreover, since $e_1e'_1=e'_1=e'_1e_1$, we have $e_1e=e=ee_1$.
Similarly, $e_2e=e=ee_2$. We get $r_1 := e'_1re'_1+ e'_1re + ere'_1
+ ere\in e_1R^{\otimes_E n}e_1=A_1^{\otimes_E n}$ and $r_2:= ere'_2
+ e'_2re + e'_2re'_2 \in e_2R^{\otimes_E n}e_2 =A_2^{\otimes_E n}$.
This shows that $g^{n}(r_1, -r_2)=r$ and $g^{n}$ is surjective.
Similar arguments show that $f^{n}$ and $g^{n}$ are $E^e$-linear.
\end{proof}

\fait{corollary}{cor cyclic}{Let $R$ be an algebra, oriented by
three idempotents $e, e'_1, e'_2$.  For any $n\geq 1$, there
exists a short exact sequence in \emph{mod}$E^e$ of the form
\begin{eqnarray*}
\qquad \SelectTips{eu}{10}\xymatrix@C=25pt{0\ar[r]&
\widetilde{C^{\otimes_{E_C} n}} \ar[r]^-{f^n} &\displaystyle
\widetilde{A_1^{\otimes_{E_1} n}}\bigoplus
\widetilde{A_2^{\otimes_{E_2} n}} \ar[r]^-{g^{n}} &
\widetilde{R^{\otimes_E n}} \ar[r] & 0}
\end{eqnarray*}}
\begin{proof}
It suffices to check that the maps $f^{n}$
and $g^{n}$ in Proposition \ref{prop suite exacte} are compatible
with the relation $\sim$ and remain respectively injective and
surjective, which is straightforward.
\end{proof}

We can now state and prove the main results of this section.

\fait{theorem}{thm cyclic homology}{Let $R$ be an algebra, oriented
by three idempotents $e, e'_1$ and $e'_2$. Moreover, let $C=eRe$,
$A_1=e_1Re_1$ and $A_1=e_2Re_2$ be as in Definition \ref{defn
oriented}. There exists a long exact sequence of cyclic homology
groups :
$$\SelectTips{eu}{10}\xymatrix@C=12pt@R=3pt{ \cdots \ar[r] & {\rm HC}_1(C) \ar[r] &
{\rm HC}_1(A_1) \oplus {\rm HC}_1(A_2) \ar[r] &
{\rm HC}_1(R) \ar[r] & \ \ \ \\ %
\ \ar[r] & {\rm HC}_0(C)\ar[r] & {\rm HC}_0(A_1) \oplus {\rm HC}_0(A_2) \ar[r] &
{\rm HC}_0(R) \ar[r] & 0}$$}
\begin{proof}
The maps $f^{n}$ and $g^{n}$ of Corollary
\ref{cor cyclic} commute with $b$, $b'$, $t$ and $N$ in
(\ref{CC_S}), so they extend to a short exact sequence of relative
bicomplexes
\begin{equation}\label{sec relative}
\SelectTips{eu}{10}\xymatrix{0\ar[r] & \widetilde{CC_{E_C}}(C) \ar[r]^-{f^\bullet} &
\widetilde{CC_{E_1}}(A_1)\oplus \widetilde{CC_{E_2}}(A_2)
\ar[r]^-{g^\bullet} & \widetilde{CC_E}(R) \ar[r] & 0}
\end{equation} The result then follows from Proposition \ref{prop quasi}.
\end{proof}

%

In Theorem \ref{thm Hochschild homology} and Theorem \ref{thm
cyclic homology}, we have obtained long exact sequences of
Hochschild homology groups and cyclic homology groups.  A deep
inspection of the natural morphisms involved in these sequences,
as well as in Connes' long exact sequence yields the following
result, whose verification is left to the reader.

\fait{theorem}{cyclic connes}{Let $R$ be an algebra, oriented by
three idempotents $e, e'_1$ and $e'_2$. Moreover, let $C=eRe$,
$A_1=e_1Re_1$ and $A_1=e_2Re_2$ be as in Definition \ref{defn
oriented}. There exists a commutative diagram $$
\SelectTips{eu}{10}\xymatrix@R=15pt@C=12pt{%
& \vdots \ar[d] & \vdots\ar[d] & \vdots\ar[d] & \vdots\ar[d]\\
\cdots \ar[r] & {\rm H}_n(R)\ar[d]\ar[r] & {\rm HC}_n(R)\ar[d] \ar[r] &
{\rm HC}_{n-2}(R)\ar[d] \ar[r] & {\rm H}_{n-1}(R)\ar[d] \ar[r] & \cdots\\
\cdots \ar[r] & {\rm H}_{n-1}(C)\ar[d]\ar[r] & {\rm HC}_{n-1}(C)\ar[d] \ar[r]
& {\rm HC}_{n-3}(C)\ar[d] \ar[r] & {\rm H}_{n-2}(C)\ar[d] \ar[r] & \cdots\\
\cdots \ar[r] & {\begin{array}{c} {\rm H}_{n-1}(A_1)\\ \oplus
\\ {\rm H}_{n-1}(A_2)\end{array}} \ar[d]\ar[r] &
{\begin{array}{c}{\rm HC}_{n-1}(A_1)\\ \oplus \\
{\rm HC}_{n-1}(A_2)\end{array}}\ar[d] \ar[r] &
{\begin{array}{c}{\rm HC}_{n-3}(A_1)\\ \oplus \\
{\rm HC}_{n-3}(A_2)\end{array}} \ar[d] \ar[r] &
{\begin{array}{c}{\rm H}_{n-2}(A_1)\\ \oplus \\
{\rm H}_{n-2}(A_2)\end{array}}\ar[d] \ar[r] &\cdots\\
\cdots \ar[r] & {\rm H}_{n-1}(R)\ar[d]\ar[r] & {\rm HC}_{n-1}(R)\ar[d] \ar[r]
& {\rm HC}_{n-3}(R)\ar[d] \ar[r] & {\rm H}_{n-2}(R)\ar[d] \ar[r] & \cdots\\
& \vdots & \vdots & \vdots & \vdots }$$ in which the rows and
columns are exact.}

We finish this section by showing that dual results hold for
cyclic cohomology groups. Let $R$ be an algebra and $S$ be a
separable subalgebra. We have seen (in our discussion before
Proposition \ref{prop quasi}) that the canonical projection
\mor{\pi}{CC(R)}{\widetilde{CC_S}(R)} is a morphism of bicomplexes
which is a quasi-isomorphism when restricted to the columns.
Because $D=\Hom_k(-,k)$ is exact, the dual morphism
\mor{D(\pi)}{D(\widetilde{CC_S}(R))}{D(CC(R)} has the same
property. Thus the map induced by $D(\pi)$ on the total complexes
is a quasi-isomorphism. We get the following.

\fait{proposition}{prop quasi coh}{Let $R$ be an algebra and $S$ be
a separable subalgebra of $R$.  Then ${\rm HC}^n(R)\simeq {\rm
H}^n(\emph{Tot}D(\widetilde{CC_S}(R)))$.}

As a consequence, to obtain a long exact sequence on cyclic homology
groups, it suffices to construct a short exact sequence of the duals
of the relative bicomplexes (for $S=E$).  But such a sequence exists
by (\ref{sec relative}) and the fact that $D$ is exact.  We have
shown the following theorem.

\fait{theorem}{thm cyclic cohomology}{Let $R$ be an algebra,
oriented by three idempotents $e, e'_1$ and $e'_2$. Moreover, let
$C=eRe$, $A_1=e_1Re_1$ and $A_1=e_2Re_2$ be as in Definition
\ref{defn oriented}. There exists a long exact sequence of cyclic
cohomology groups:
$$\SelectTips{eu}{10}\xymatrix@C=12pt@R=3pt{0 \ar[r] & {\rm HC}^0(R) \ar[r] & {\rm HC}^0(A_1)
\oplus {\rm HC}^0(A_2) \ar[r] &
{\rm HC}^0(C) \ar[r] & \ \ \ \\ %
\ \ar[r] & {\rm HC}^1(R)\ar[r] & {\rm HC}^1(A_1) \oplus {\rm
HC}^1(A_2) \ar[r] & {\rm HC}^1(C) \ar[r] & \cdots}$$}

Here again, the long exact sequences obtained in Theorem \ref{thm
Hochschild cohomology} and Theorem \ref{thm cyclic cohomology}
embed with the cohomological Connes' exact sequence, so that we
get the following result, whose verification is left to the
reader.
\fait{theorem}{cyclic connes coh}{Let $R$ be an algebra, oriented by
three idempotents $e, e'_1$ and $e'_2$. Moreover, let $C=eRe$,
$A_1=e_1Re_1$ and $A_1=e_2Re_2$ be as in Definition \ref{defn
oriented}. There exists a commutative diagram $$
\SelectTips{eu}{10}\xymatrix@R=15pt@C=9pt{%
& \vdots \ar[d] & \vdots\ar[d] & \vdots\ar[d] & \vdots\ar[d]\\
\cdots \ar[r] & \mbox{\footnotesize{${\rm H}^n(R, D(R))$}}\ar[d]\ar[r] &
\mbox{\footnotesize{${\rm HC}^{n-1}(R)$}}\ar[d] \ar[r] &
\mbox{\footnotesize{${\rm HC}^{n+1}(R)$}}\ar[d] \ar[r] &
\mbox{\footnotesize{${\rm H}^{n+1}(R, D(R))$}}\ar[d] \ar[r] & \cdots\\
\cdots \ar[r] & \mbox{\footnotesize{${\begin{array}{c} {\rm H}^{n}(A_1,
D(A_1))\\ \oplus
\\ {\rm H}^{n}(A_2, D(A_2))\end{array}}$}} \ar[d]\ar[r] &
\mbox{\footnotesize{${\begin{array}{c}{\rm HC}^{n-1}(A_1)\\ \oplus \\
{\rm HC}^{n-1}(A_2)\end{array}}$}}\ar[d] \ar[r] &
\mbox{\footnotesize{${\begin{array}{c}{\rm HC}^{n+1}(A_1)\\ \oplus \\
{\rm HC}^{n+1}(A_2)\end{array}}$}} \ar[d] \ar[r] &
\mbox{\footnotesize{${\begin{array}{c}{\rm H}^{n+1}(A_1, D(A_1))\\
\oplus \\ {\rm H}^{n+1}(A_2, D(A_2))\end{array}}$}}\ar[d] \ar[r]
&\cdots\\
\cdots \ar[r] & \mbox{\footnotesize{${\rm H}^{n}(C,D(C))$}}\ar[d]\ar[r]
& \mbox{\footnotesize{${\rm HC}^{n-1}(C)$}}\ar[d] \ar[r] &
\mbox{\footnotesize{${\rm HC}^{n+1}(C)$}}\ar[d] \ar[r] &
\mbox{\footnotesize{${\rm H}^{n+1}(C,D(C))$}}\ar[d] \ar[r] & \cdots\\
\cdots \ar[r] &
\mbox{\footnotesize{${\rm H}^{n+1}(R,D(R))$}}\ar[d]\ar[r] &
\mbox{\footnotesize{${\rm HC}^{n}(R)$}}\ar[d] \ar[r] &
\mbox{\footnotesize{${\rm HC}^{n+2}(R)$}}\ar[d] \ar[r] &
\mbox{\footnotesize{${\rm H}^{n+2}(R, D(R))$}}\ar[d] \ar[r] & \cdots\\
& \vdots & \vdots & \vdots & \vdots }$$ in which the rows and
columns are exact.}

%
%
\section{Simplicial (co)homology and fundamental group}
    \label{Simplicial and fundamental}

In this section we consider algebras of the form $R=kQ_R/I_R$. More
precisely we study the simplicial (co)homology groups of $R$ (in
case they admit semi-normed basis) as well as the fundamental group
$\pi_1(Q_R,I_R)$. As before, we show the existence of Mayer-Vietoris
long exact sequence relating the simplicial (co)homology groups.
Then we prove that under some additional assumptions, it is possible
to express the fundamental group of $(Q_R,I_R)$ in terms of those of
$(Q_C,I_C)$, $(Q_{A_1},I_{A_1})$ and $(Q_{A_2},I_{A_2})$.  We
further apply this result to the computation of the first Hochschild
cohomology group of some schurian oriented algebras.

We stress here that the validity of one of the conditions (\ref{case
1})-(\ref{case 3}) is not necessary in the sequel, and so the
results contained in this section apply for a slightly more general
family of algebras, that is those satisfying (possibly only)
condition (\ref{case 0}), for instance the algebra of Examples
\ref{ex tournicoti} (b).  Nevertheless, we will not make any
distinction and we will keep going with the expression "oriented
algebras".
Also, all algebras considered in this section are bound quiver
algebras, that is of the form $kQ/I$ for some bound quiver $(Q,I)$
(see Sections \ref{Algebras and quivers} and \ref{The case where C
is a core algebra} for more details).

\subsection{Simplicial (co)homology}
    \label{simpA}

Given an algebra $R=kQ/I$ and a pair $(x, y)$ of vertices of $Q$,
let $_y\mathbb B_x$ be a basis of the $k$-space $\varepsilon_xR\varepsilon_y$  and
$\mathbb B = \bigcup_{(x,y) \in Q_0 \times Q_0}{} _y\mathbb B_x$.
We say that $\mathbb B$ is a \textbf{semi-normed basis}
\cite{MdlP99} if:
\begin{enumerate}
\item[(a)] $\varepsilon_x \in {} _x\mathbb B_x$ for every vertex $x \in
Q_0$; %
\item[(b)] $\alpha + I  \in {}_y\mathbb B_x$ for every arrow
\mor{\alpha}{x}{y};%
\item[(c)] For $\sigma$ and $\sigma '$ in $\mathbb
B$, either $\sigma\sigma ' = 0$ or there exists $0 \neq
\lambda_{\sigma,\sigma '} \in k$, and $b(\sigma,\sigma ') \in
\mathbb B$ such that $\sigma\sigma ' =\lambda_{\sigma,\sigma
'}b(\sigma,\sigma ')$.
\end{enumerate}

Assume $R=kQ_R/I_R$ admits a semi-normed basis.  Following
\cite{MVP83, MdlP99}, one can associate to $R$ the chain complex
$({\rm SC}_{\bullet}(R), d)$ as follows: ${\rm SC}_0(R)$ and ${\rm
SC}_1(R)$ are the free abelian groups with basis $Q_0$ and
$\mathbb B$, respectively. For $n \geq 2$, let ${\rm SC}_n(R)$ be
the free abelian group with basis the set of $n$-tuples
$(\sigma_1, \sigma_2, \ldots, \sigma_n)$ of $\mathbb{B}^n$ such
that $\sigma_1 \sigma_2 \cdots \sigma_n \neq 0$, and $\sigma _i
\neq e_j$ for all $i, j \in Q_0$. The differential \mor{d_n}{{\rm
SC}_n(R)}{{\rm SC}_{n-1}(R)} is defined on the basis elements by the rules $d_1(\sigma)
= y-x $ for $\sigma \in \ _y\mathbb B_x$  and
\begin{eqnarray*}
d_n(\sigma_1, \sigma_2, \dots, \sigma_n)&=& (\sigma_2, \ldots,
\sigma_n)\\
& + & \sum_{j=1}^{n-1}(-1)^j(\sigma_1, \ldots, b(\sigma_j,
\sigma_{j+1}), \ldots, \sigma_n)\\
&+&(-1)^n (\sigma_1, \ldots, \sigma_{n-1})
\end{eqnarray*}
for $n \geq 1$ and $(\sigma_1, \sigma_2, \dots, \sigma_n)\in {\rm
SC}_n(R)$.

The $n$th \textbf{simplicial homology group} of $R$ (with respect
to the basis $\mathbb{B}$) is the $n$th homology group of the
chain complex $({\rm SC}_\bullet(R),d_\bullet)$, and we denote it by
$\SH_n(R)$. Given an abelian group $G$, the
$n$th-\textbf{simplicial cohomology group} $\SH^n(R, G)$ of $R$
with coefficients in $G$ is the $n$th cohomology group of the
complex obtained by applying the functor $\Hom_{\mathbb Z} (-, G)$
on the chain complex $({\rm SC}_\bullet(R),d_\bullet)$.

Observe that the above chain complexes depend essentially on the
way $\mathbb{B}$ is related to $I$. Hence, as for fundamental
groups (see Section \ref{simpB}), different presentations of $R$
may lead to different simplicial (co)homology groups. See \cite[(5.2)]{Bus04} for an example.

If $R$ is oriented by $e, e'_1$ and $e'_2$ with $\mathbb B_1$ and
$\mathbb B_2$ semi-normed bases of $A_1$ and $A_2$, respectively,
then $\mathbb{B}_R = \mathbb B_1 \cup \mathbb B_2$ and
$\mathbb{B}_C=\mathbb{B}_1\cap\mathbb{B}_2$ are semi-normed bases
of $R$ and $C$, respectively.  This leads to the following
results.

\fait{proposition}{SimpHom}{Let $R$ be oriented by $e, e'_1$ and
$e'_2$. Moreover, assume that $A_1$ and $A_2$ admit semi-normed
basis. Then there exists a long exact sequence of simplicial
homology groups: $$\SelectTips{eu}{10}\xymatrix@C=12pt@R=12pt{ \cdots \ar[r] &
\SH_n(C) \ar[r] & \SH_n(A_1) \oplus \SH_n(A_2) \ar[r] & \SH_n(R)
\ar[r] & \SH_{n-1}(C)\ar[r] & \cdots}$$ }
\begin{proof}
By the above discussion, $R$ and $C$
admit semi-normed bases. For $i\geq 0$, define %
$$p_i : \SelectTips{eu}{10}\xymatrix@R=0pt{\qquad \quad {\rm SC}_i(A_1) \oplus {\rm
SC}_i(A_2) \ar[r] & {\rm SC}_i(R) \qquad \qquad \qquad \qquad
\\
((\beta_1, \dots, \beta_i),(\gamma_1, \dots, \gamma_i))
\ar@{|->}[r] &  (\beta_1, \dots, \beta_i)- (\gamma_1, \dots,
\gamma_i)}$$
The map $p_\bullet$ is a morphism of complexes. Moreover, let
$\alpha = (\alpha_1, \dots,\alpha_i)$ be a generator of
$\mbox{SC}_i(R)$. Since $e'_1Re'_2= e'_2Re'_1=0$, then, for each
$j$, we have $\alpha_j\in \mathbb{B}_1$ or $\alpha_j\in
\mathbb{B}_2$. Consequently, $\alpha \in {\rm SC}_i(A_1)$ or $\alpha
\in {\rm SC}_i(A_2)$, whence $p_i$ is surjective. Moreover, the
kernel of $p_i$ is exactly ${\rm SC}_i(C)$. We obtain a short exact
sequence of complexes
 \begin{equation}  \label{JC}
\SelectTips{eu}{10}\xymatrix@C=12pt@R=12pt{ 0 \ar[r]  & {\rm SC}_\bullet(C) \ar[r]
& {\rm SC}_\bullet(A_1)\oplus {\rm SC}_\bullet(A_2) \ar[r] & {\rm
SC}_\bullet(R) \ar[r] & 0}
\end{equation}
%
%
from which we deduce the result.
\end{proof}

\fait{proposition}{SimpCohom}{Let $R$ be oriented by $e, e'_1$ and
$e'_2$. Moreover, assume that $A_1$ and $A_2$ admit semi-normed
basis. Then, for each abelian group $G$, there exists a long exact
sequence of simplicial cohomology groups:
$$\SelectTips{eu}{10}\xymatrix@C=7pt@R=12pt{ \cdots \ar[r]  & \SH^n(R, G) \ar[r] &
\SH^n(A_1, G) \oplus \SH^n(A_2, G) \ar[r] & \SH^n(C, G) \ar[r] &
\SH^{n+1}(R, G)\ar[r] & \cdots}$$}
\begin{proof}
Since ${\rm SC}_i(R)$ is a free abelian
group for all $i$, the sequence (\ref{JC}) splits  and so the short
sequence of complexes obtained by applying $\mbox{Hom}_\mathbb{Z}(-,
G)$  to the sequence (\ref{JC}) is exact.
\end{proof}

\subsection{Fundamental group}
    \label{simpB}

Let $(Q,I)$ be a bound quiver. Recall from \cite{MVP83} that, given
two vertices $x, y \in Q_0$, a relation $\rho = \sum_{i=1}^{n}
\lambda_i\omega_i \in I\cap kQ(x,y)$ is called \textbf{minimal}
provided $n \geq 2$ and $\sum_{i \in J} \lambda_i\omega_i \notin
I\cap kQ(x,y)$ for every  proper subset $J$ of $\{1, 2, \dots, n\}$.
Let $Q'_1$ be the union of $Q_1$ with the set of inverse arrows
\mor{\alpha^{-1}}{y}{x} for each arrow \mor{\alpha}{x}{y} in $Q_1$,
a \textbf{walk} in $(Q_0, Q_1)$ is a path in $(Q_0, Q'_1)$. The
\textbf{homotopy relation} in the set of walks on $Q$ is the
smallest equivalence relation such that:
\begin{enumerate}
\item[(a)] For each arrow \mor{\alpha}{x}{y} $\in Q_1$, $\alpha \alpha^{-1} \sim
\epsilon_x$ and $\alpha^{-1}\alpha \sim \epsilon_y$;
\item[(b)] If $\sum_{i=1}^{n} \lambda_i\omega_i$ is a minimal relation, then
$\omega_j \sim \omega_k$ for all $1 \leq j,k\leq n$;
\item[(c)] If $u,v,w, w'$ are walks such that $w \sim w'$, then $uwv \sim uw'v$,
whenever these compositions are defined.
\end{enumerate}

A minimal relation $\rho$ is \textbf{fundamental} if there does not
exist a nontrivial path $u$ and a minimal relation $\omega$ such
that $\rho=\omega u$ or $\rho=u\omega$. One can see that the above
homotopy relation does not change if we replace the condition (b) by
\begin{enumerate}
\item[(b')] If $\sum_{i=1}^{n} \lambda_i\omega_i$ is a fundamental minimal relation, then
$\omega_j \sim \omega_k$ for all $1 \leq j,k\leq n$.
\end{enumerate}

When dealing with homotopy, it is thus sufficient to consider only
the fundamental minimal relations. Here, we give examples of
fundamental minimal relations.

\fait{example}{ex fund}{\emph{%
Let $Q$ be the quiver
$$\SelectTips{eu}{10}\xymatrix@C=15pt@R=8pt{ & 2\ar[dl]_{\alpha_1} & & & 6\ar[dl]_{\delta_1}\\%
1 & & 4 \ar[ul]_{\beta_1} \ar[dl]^{\beta_2} & 5 \ar[l]_\gamma & &
8 \ar[ul]_{\varepsilon_1} \ar[dl]^{\varepsilon_2} \\%
& 3 \ar[ul]^{\alpha_2} & & & 7 \ar[ul]^{\delta_2}}$$
bound by the relations $\beta_1 \alpha_1 - \beta_2\alpha_2 =
\varepsilon_1 \delta_1 - \varepsilon_2 \delta_2
=\delta_1\gamma\beta_1=\delta_2\gamma\beta_2=0$.  The minimal
fundamental relations are $\beta_1\alpha_1-\beta_2\alpha_2$ and
$\varepsilon_1 \delta_1 - \varepsilon_2 \delta_2$, while the
relations $\gamma(\beta_1\alpha_1-\beta_2\alpha_2)$ and
$(\varepsilon_1 \delta_1 - \varepsilon_2 \delta_2)\gamma$ are
minimal but not fundamental.}}

For a fixed vertex $x\in Q_0$, let $\pi_1(Q, x)$ be the fundamental
group of the underlying graph of $Q$ at the vertex $x$, that is the
set of homotopy classes of closed walks on $Q$ starting and ending
at $x$. It is well-known (see \cite{R88}, for instance) that
$\pi_1(Q, x)$ is isomorphic to the free multiplicative nonabelian
group on $\chi (Q) = |Q_1|-|Q_0|+1$ generators. Now, let $N(Q, I,
x)$ be the normal subgroup of $\pi_1(Q, x)$ generated by the
elements of the form $wuv^{-1}w^{-1}$ where $w$ is a path starting
at $x$ and $u,v$ are two homotopic paths. The \textbf{fundamental
group at $x$} of the bound quiver $(Q, I)$ is defined to be
$\pi_1(Q,I,x) = \pi_1(Q)/ N(Q,I,x)$. If $Q$ is connected, this
definition does not depend on the choice of $x$. Accordingly, we
write $\pi_1(Q,I)=\pi_1(Q)/N(Q,I)$, and call it the
\textbf{fundamental group} of the bound quiver $(Q,I)$. Similarly,
we write $\pi_1(Q)$ instead of $\pi_1(Q, x)$.

On other hand, following \cite{Bus04}, to a bound quiver $(Q,I)$, we
associate a CW complex $\mathcal B(Q,I)$ (called its
\textbf{classifying space}) in the following way: the $0$-cells are
given by $Q_0$ and the $n$-cells are given by $n$-tuples
$(\overline{\sigma_1}, \dots, \overline{\sigma_n})$ of homotopy
classes of nontrivial paths of $(Q ,I)$ such that the composition
$\sigma_1\sigma_2\cdots\sigma_n$ is a path of $(Q,I)$ not in $I$.
Then, by \cite[(3.3)]{Bus04}, there is an isomorphism of groups
\begin{equation}\label{iso pi1}
 \pi_1(Q,I) \simeq  \pi_1({\mathcal B}(Q,I))
\end{equation}

%
We shall use this fact, together with Van Kampen's theorem, in the following proof.

\fait{proposition}{VK}{Let $R=kQ_R/I_R$ be a bound quiver algebra,
oriented by three idempotents $e, e'_1$ and $e'_2$. Moreover, assume
that $C=kQ_{C}/I_C$, $A_1=kQ_{A_1}/I_1$ and $A_2=kQ_{A_2}/I_2$ and
that, for each fundamental minimal relation $\omega$ on $(Q_R,
I_R)$, all the paths occurring in $\omega$ are either contained in
$Q_C$ or contain no arrows from $Q_C$. Let $Q_C=\coprod_{j=1}^m
Q_{C_j}$ with each $Q_{C_j}$ connected, and set $I_{C_j}=I_R\cap
kQ_{C_j}$. Then,
$$\pi_1(Q_R,I_R) \simeq
\frac{\pi_1(Q_{A_1},I_{A_1})\coprod\pi_1(Q_{A_2},I_{A_2})
\coprod L_{m-1}}{\coprod_{j=1}^m\pi_1(Q_{C_j},I_{C_j})}$$
where $\coprod$ denotes the free product of groups and $L_{m-1}$
is the free group in $m-1$
generators.}
\begin{proof}
Let $n_1, n_2, n_R, n$ and $n_j$ be the number of vertices in
$Q_{A_1}, Q_{A_2}, Q_R, Q_C$ and $Q_{C_j}$ respectively. In
particular, $n_R = n_1 + n_2 -n$. Let $T_{C_j}$ be a maximal tree in
$Q_{C_j}$, $T_{C}$ be the disjoint union of the $T_{C_j}$ and for $i
\in \{1,2\}$, let $T_i$ be a maximal tree in $Q_{A_i}$ such that
$T_1 \cap T_2 = T_C$. Such maximal trees always exist.  In addition,
$|{(T_C)}_1|=\sum_{j=1}^m |{(T_{C_j})}_1|=\sum_{j=1}^m(n_j-1)=n-m$.
For $i=1,2$, set
$$\hat{Q}_{A_i}=Q_{A_i} \coprod_{T_C}T_k \textrm{ with } k\neq i
\textrm{ and } \hat{I}_{i}=I_R \cap k\hat{Q}_{A_i}.$$
Observe that $Q_R=\hat{Q}_{A_1} \coprod_{Q_C} \hat{Q}_{A_2}$.  In
addition, $|{(\hat{Q}_{A_1})}_0|=n_R$ and
$$ |{(\hat{Q}_{A_1})}_1|= |{(Q_{A_1})}_1| +
|{(T_2)}_1|-|{(T_C)}_1| = |{(Q_{A_1})}_1| + (n_2-1)-(n-m)$$
and a similar expression holds for $|{(\hat{Q}_{A_2})}_1|$. So,
$\pi_1(\hat{Q}_{A_i})$ is the free group in $\chi(\hat{Q}_{A_i})$
generators, where $ \chi(\hat{Q}_{A_i}) = \chi(Q_{A_i}) +(m-1)$.
Then, $\pi_1(\hat{Q}_{A_i}) \simeq \pi_1(Q_{A_i}) \coprod L_{m-1}$.

Moreover, by our assumption on the fundamental minimal relations,
we have an inclusion of groups
\mor{\iota_i}{\pi_1(\hat{Q}_{A_i},\hat{I}_{i})}{\pi_1(Q_R,I_R)}.
Also, the restriction to $(Q_{A_i},I_{i})$ of the homotopy
relation on $(\hat{Q}_{A_i},\hat{I}_{i})$ coincides with the
homotopy relation on $(Q_{A_i},I_{i})$. Thus $\hat{I}_{i} = I_{i}
+ J_{i}$ for some monomial ideal $J_{i}$ and
\begin{eqnarray*}
\pi_1(\hat{Q}_{A_i}, \hat{I}_{i})& \simeq &
\pi_1(\hat{Q}_{A_i})/N(\hat{Q}_{A_i}, \hat{I}_{i})
\\& \simeq & (\pi_1(Q_{A_i}) \coprod L_{m-1})/N(\hat{Q}_{A_i},
\hat{I}_{i}) \\& \simeq &[\pi_1(Q_{A_i}) /N(Q_{A_i}, I_{i})]
\coprod L_{m-1}
\\ & \simeq & \pi_1(Q_{A_i}, I_{i}) \coprod L_{m-1}.
\end{eqnarray*}
Now set $Q= \hat{Q}_{A_1} \cap \hat{Q}_{A_2}$ and $I=kQ \cap I_R$.
We have
\begin{eqnarray*}
\chi (Q)& = & (|{(T_1)}_1|+|{(T_2)}_1|-2|{(T_C)}_1|+|{(Q_C)}_1|)-n_R+1%
\\& = & |{(Q_C)}_1|-(n-m)+(m-1)
\\& = & (\sum_{j=1}^m\chi({Q_C}_j))+(m-1).
\end{eqnarray*}
and so $\pi_1(Q)=(\coprod_{j=1}^m \pi_1(Q_{C_j}))\coprod L_{m-1}$.
Moreover, by our assumption on the fundamental minimal relations, we
have an inclusion of groups
\mor{\kappa_i}{\pi_1(Q,I)}{\pi_1(\hat{Q}_{A_i},\hat{I}_{i})} and $I
= I_{C} + J_I$ for some monomial ideal $J_I$. So, $\pi_1(Q,I)=
L_{m-1} \coprod $ \newline $(\coprod_{j=1}^m \pi_1(Q_{C_j},
I_{C_j}))$.

The next step is to use Van Kampen's theorem for bound quivers. For
the sake of this, let $\hat{\mathcal B_1}$, $\hat{\mathcal B_2}$,
$\mathcal B_R$ and $\mathcal B$ be the classifying spaces of
$(\hat{Q}_{A_1}, \hat{I}_{1})$, $(\hat{Q}_{A_2}, \hat{I}_{2})$,
$(Q_R, I_R)$, $(Q, I)$, respectively. Then, $\mathcal B_R
=\hat{\mathcal B_1} \cup \hat{\mathcal B_2}$ and $\mathcal
B=\hat{\mathcal B_1} \cap \hat{\mathcal B_2}$, where $\mathcal{B}$
is connected. Applying Van Kampen's theorem (see \cite{R88}) and the
isomorphism of fundamental groups (\ref{iso pi1}), we have
\begin{eqnarray*}
\pi_1(Q_R,I_R) & \simeq & \pi_1(\mathcal{B}_R) \\
& \simeq &
\frac{\pi_1(\mathcal{B}_1)\coprod\pi_1(\mathcal{B}_2)}{\pi_1(\mathcal{B})}\\
& \simeq & \frac{\pi_1(\hat{Q}_{A_1},\hat{I}_{1})\coprod
\pi_1(\hat{Q}_{A_2},\hat{I}_{2})}{\pi_1(Q,I)}\\
& \simeq & \frac{(\pi_1(Q_{A_1},I_{1})\coprod L_{m-1})\coprod
(\pi_1(Q_{A_2},I_{2}) \coprod
L_{m-1})}{(\coprod_{j=1}^m\pi_1(Q_{C_j},I_{C_j})) \coprod
L_{m-1}}\\
& \simeq &
\frac{\pi_1(Q_{A_1},I_{1})\coprod\pi_1(Q_{A_2},I_{2})\coprod
L_{m-1}}{\coprod_{j=1}^m\pi_1(Q_{C_j},I_{C_j})}
\end{eqnarray*}
\end{proof}

We note that the main tool used in this proof is that the
restriction of the homotopy relation on $(Q_R,I_R)$ to $(Q_C,I_C)$
coincides with the homotopy relation on $(Q_C,I_C)$. This proof
may thus clearly be generalized to this more general context.

The above formula allows to compute the fundamental group of
several oriented algebras, but it does not include all oriented
algebras. The following examples illustrate this point.

\fait{examples}{ex pi}{\emph{%
\begin{enumerate}
\item[(a)] Let $R$ be the algebra given by the bound quiver of Example \ref{ex fund}.
Let $e$, $e'_1$ and $e'_2$ be the sums of the primitive
idempotents corresponding to the sets of vertices $\{4,5\}$,
$\{1,2,3\}$ and $\{6,7,8\}$ respectively. Clearly $R$ is oriented
by $e, e'_1$ and $e'_2$.  Moreover, $R$ satisfies to the
conditions of Proposition \ref{VK}. Since straightforward
computations show that $\pi_1(Q_{C},I_C)$, $\pi_1(Q_{A_1},I_1)$
and $\pi_1(Q_{A_2}, I_2)$ are trivial, so is $\pi_1(Q_R,I_R)$ by
Proposition \ref{VK}.
\item[(b)] Let $R$ be the algebra given by the quiver
$$\SelectTips{eu}{10}\xymatrix@C=15pt@R=10pt{1 \ar[r]^{\alpha_1}
& 3 \ar@<.5ex>[d]^{\gamma} \ar@<-.5ex>[d]_{\beta} & 2\ar[l]_{\alpha_2}\\%
& 4}$$
bound by the relations $\alpha_i(\beta-\gamma)$, for $i=1,2$. Let
$e$, $e'_1$ and $e'_2$ be the sums of the primitive idempotents
corresponding to the sets of vertices $\{3,4\}$, $\{1\}$ and
$\{2\}$ respectively. Clearly $R$ is oriented by $e, e'_1$ and
$e'_2$. However, $R$ does not satisfy the conditions of
Proposition \ref{VK} because the fundamental minimal relations
$\alpha_i(\beta-\gamma)$ contain arrows from $A'_i$ and $C$.
Moreover, $\pi_1(Q_R,I_R)$, $\pi_1(Q_{A_1},I_1)$ and
$\pi_1(Q_{A_2},I_2)$ are trivial, while $\pi_1(Q_C,I_C)\simeq
\mathbb{Z}$. So the formula of Proposition \ref{VK} does not hold.
\end{enumerate}}}

We end this section with an application of the above result to the
computation of the first Hochschild cohomology group for some
schurian oriented algebras (see Examples \ref{core algebras} (a)).
Recall that if $A=kQ/I$ is a schurian algebra, then
${\rm H}^1(A)\simeq \mbox{Hom}(\pi_1(Q, I), k^+)$ by \cite{PS01,
AdlP96}, where $k^+$ is the (abelian) additive group of $k$ and
the homomorphisms are group homomorphisms. This leads to the
following result:

\fait{proposition}{prop H1}{Assume $R$ is oriented by $e, e'_1$
and $e'_2$ and satisfies to the conditions of Proposition
\ref{VK}. If $R$ is schurian and
$\coprod_{i=1}^m\pi_1(Q_{C_i},I_{C_i})$ is trivial, then
${\rm H}^1(R)\simeq {\rm H}^1(A_1)\oplus {\rm H}^1(A_2) \oplus k^{m-1}$.}
\begin{proof}
Since $R$ is schurian, then so are $A_1$,
$A_2$ and $C$.  Therefore%
$$\begin{array}{rcl}
  {\rm H}^1(R) & \simeq & \mbox{Hom}(\pi_1(Q_R,I_R), k^+) \\
   & \simeq &  \mbox{Hom}(\pi_1(Q_{A_1},I_{A_1})\amalg \pi_1(Q_{A_2},I_{A_2}) \amalg L_{m-1}, k^+) \\
   & \simeq &  \mbox{Hom}(\pi_1(Q_{A_1}, I_1), k^+) \oplus \mbox{Hom}(\pi_1(Q_{A_2}, I_2), k^+)
   \oplus k^{m-1}\\
   & \simeq & {\rm H}^1(A_1)\oplus {\rm H}^1(A_2) \oplus k^{m-1}
\end{array} $$
\end{proof}

\fait{example}{ex tournicoti 2}{\emph{%
Let $R$ be the algebra of Examples \ref{ex tournicoti} (b). Theorem
\ref{thm Hochschild cohomology} cannot be used to compute the
Hochschild cohomology groups of $R$.  However, $R$ is schurian and
clearly satisfies (\ref{case 0}) and the conditions of Proposition
\ref{VK} since $C$ is semisimple. Thus ${\rm H}^1(R)\simeq {\rm
H}^1(A_1)\oplus
{\rm H}^1(A_2)\oplus k \simeq k$ by Proposition \ref{prop H1}. %
}}

\medskip

%
%
ACKNOWLEDGEMENTS. This paper was started while the first author
visited the Universities of Sherbrooke and Bishop's under a
post-doctoral fellowship, and while the second and the third author
respectively had an ISM and NSERC graduate studies fellowship;
unfortunately, the write-up of the results has been delayed for
quite a while. The authors warmly thank Ibrahim Assem for his
precious help and patience.

\medskip

%
%
\providecommand{\bysame}{\leavevmode\hbox
to3em{\hrulefill}\thinspace}
\providecommand{\MR}{\relax\ifhmode\unskip\space\fi MR }
\providecommand{\MRhref}[2]{%
  \href{http://www.ams.org/mathscinet-getitem?mr=#1}{#2}
} \providecommand{\href}[2]{#2}

\end{document}